\documentclass[12pt,journal,letterpaper,onecolumn]{IEEEtran}

\usepackage{graphicx}
\graphicspath{{pics/}{figs/}}
\usepackage{graphicx}
\usepackage{subfigure}
\usepackage{amssymb}
\usepackage{amsmath}

\usepackage{graphicx}
\graphicspath{{pics/}{figs/}}
\usepackage{graphicx}
\usepackage{subfigure}
\usepackage{amssymb}
\usepackage{amsmath}

\newtheorem{defn}{\noindent $\mathbf{Definition}$}[section]

\newtheorem{thm}[defn]{$\noindent \mathbf{Theorem}$}

\begin{document}
\title{Computing Quasiconformal Maps on Riemann surfaces using Discrete Curvature Flow}
\author{W. Zeng \thanks{W. Zeng and L.M. Lui contributed equally for this work}, L.M. Lui \thanks{}, F. Luo,  J.S. Liu, T.F. Chan, S.T. Yau, X.F. Gu}%$^*$Wei Zeng$^\dag$ \thanks{$^\dag$Department of Computer Science, State University of New York at Stony Brook \{zengwei@cs.stonybrook.edu\}}, $^*$Lok Ming Lui$^{\ddag}$ \thanks{$\ddag$Department of Mathematics, Harvard University and UCLA\{malmlui@math.harvard.edu\}}, Feng Luo$^{\S}$, Tony Chan$^{\natural}$, Shing-Tung Yau$^{\P}$, Xianfeng Gu$^{\flat}$}

\markboth{Preprint}{Shell
\MakeLowercase{\textit{et al.}}: Beltrami Equation on Riemann surfaces}

% If you want to put a publisher's ID mark on the page
% (can leave text blank if you just want to see how the
% text height on the first page will be reduced by IEEE)
%\pubid{0000--0000/00\$00.00~\copyright~2002 IEEE}

%\specialpapernotice{(Invited Paper)}

\maketitle
\begin{abstract}
Surface mapping plays an important role in geometric processing. They induce both area and angular distortions. If the angular distortion is bounded, the mapping is called a {\it quasi-conformal} map. Many surface maps in our physical world are quasi-conformal. The angular distortion of a quasi-conformal map can be represented by Beltrami differentials. According to quasi-conformal Teichm\"uller theory, there is an 1-1 correspondence between the set of Beltrami differentials and the set of quasi-conformal surface maps. Therefore, every quasi-conformal surface map can be fully determined by the Beltrami differential and can be reconstructed by solving the so-called Beltrami equation.

In this work, we propose an effective method to solve the Beltrami equation on general Riemann surfaces. The solution is a quasi-conformal map associated with the prescribed Beltrami differential. We firstly formulate a discrete analog of quasi-conformal maps on triangular meshes. Then, we propose an algorithm to compute discrete quasi-conformal maps. The main strategy is to define a discrete auxiliary metric of the source surface, such that the original quasi-conformal map becomes conformal under the newly defined discrete metric. The associated map can then be obtained by using the discrete Yamabe flow method. Numerically, the discrete quasi-conformal map converges to the continuous real solution as the mesh size approaches to 0. We tested our algorithm on surfaces scanned from real life with different topologies. Experimental results demonstrate the generality and accuracy of our auxiliary metric method.
\end{abstract}


\section{Introduction}
Mapping between surfaces plays a fundamental role in digital geometric processing. In general, surface mappings introduce distortions, which can be classified as area distortion and angular distortion. Mappings without angular distortions are called conformal mappings. Last several years,
there has been fast development of various techniques for computing conformal mappings and their applications in geometric
processing. However, conformal mappings are not common in practice. Many mappings in our physical world are quasi-conformal, which introduce bounded angular distortion. For example, deformations of elastic shapes are quasi-conformal, such as human expression change, deformations of human organs, etc. In order to model surface mappings in the real world more effectively, it is crucial to study quasi-conformal mappings which allow for a much wider domain of applications.

\begin{figure}
\centering
\includegraphics[height=3.8in]{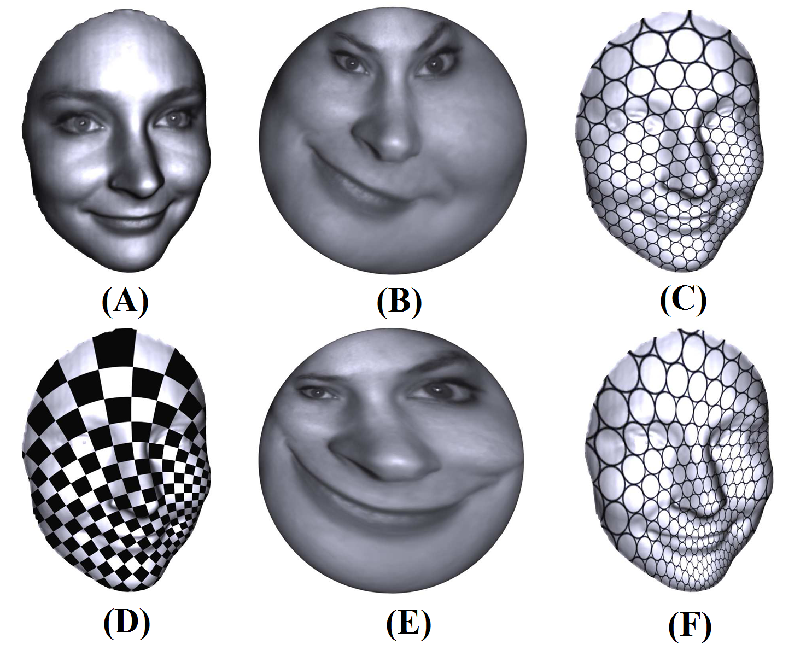}
\vspace{-2mm} \caption{Conformal and Quasi-conformal maps for a
topological disk. (A) is the original face. (B) is the conformal mapping of (A). (C) is the circle packing induced by (B). (D) is the checkerboard texture induced by (B). (E) is the quasi-conformal mapping. (F) is the circle packing induced by (E). \label{fig:riemann_mapping}} \vspace{-2mm}
\end{figure}

\begin{figure}
\centering
\includegraphics[width=2.5in]{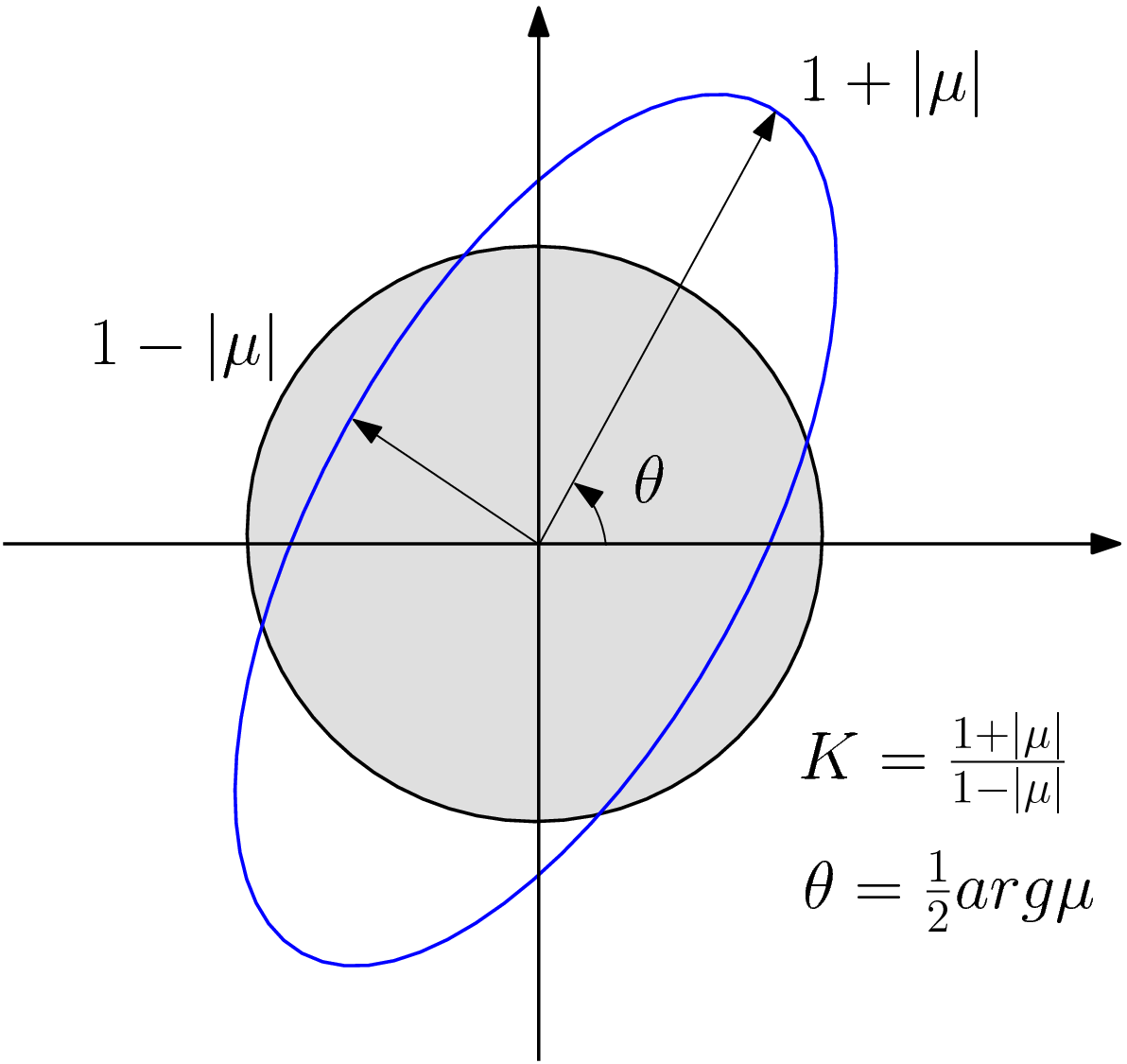}
\vspace{-2mm} \caption{Illustration of how the Beltrami coefficient $\mu$ measures the distortion by a quasi-conformal mapping that maps
a small circle to an ellipse with dilation $K$.
\label{beltrami}}
\end{figure}

The theory of quasi-conformal mappings is nearly 80 years old and has been firstly studied by Ahlfors \cite{Ahlfors}, Grotzch \cite{Grotzsch}, Morrey \cite{Morrey} and Lavrentjev \cite{Lavrentjev}. Quasi-conformal mappings can be viewed as a generalization of conformal mappings. Figure \ref{fig:riemann_mapping} illustrates the difference between a conformal map and a quasi-conformal map. Angular distortion can be characterized in the following intuitive way. Geometrically, a conformal mapping maps an infinitesimal circle on the source surface to a circle on the target surface, as shown in the first row. A quasi-conformal mapping maps an infinitesimal circle to an ellipse, as shown in frame (E) and (F). The eccentricity and the orientation of the ellipse can be represented by a complex valued function, the so-called \emph{Beltrami coefficient $\mu$}. Specifically, the ratio between the two axes of the ellipse is given by $\frac{1+|\mu(z)|}{1-|\mu(z)|}$, and the orientation of the axis is related to $arg \mu(z)$ (see figure \ref{beltrami}).

Beltrami coefficient is defined on a local chart. Globally, Beltrami coefficient is represented by the \emph{Beltrami
differential}, which is independent of the choice of local parameters. According to quasi-conformal Teichm\"uller theory, for general surfaces, there is an one to one correspondence between the set of quasi-conformal maps and the set of Beltrami differentials. In other words, every quasi-conformal map can be fully determined by the Beltrami differentials and is unique up to a finite dimensional group. Conversely, given a particular Beltrami differential $\mu(z) \frac{d\overline{z}}{dz}$, we can reconstruct the quasi-conformal maps associated to $\mu(z) \frac{d\overline{z}}{dz}$. The Beltrami differential captures the most essential information of the surface mappings. Therefore, by adjusting $\mu(z) \frac{d\overline{z}}{dz}$, we can reconstruct a surface mapping with desired properties.

Quasi-conformal mappings have been studied extensively in complex analysis \cite{Ahlfors,Grotzsch,Ahlfors2,Lehto}. Applications can be found in different areas such as differential equations, topology, Riemann mappings, complex dynamics, grid generation and so on \cite{Bers, Gamelin, Belinskii, Bers2, Bers3, Bers4, Daripa, Mastin}. Despite the rapid development in the theory of quasi-conformal geometry, the progress on computing quasi-conformal mappings numerically has been very slow. In fact, developing an effective numerical algorithm to compute quasi-conformal mapping remains a challenging problem.

Recently, there has been a few work on numerical quasi-conformal mapping techniques on the complex plane based on solving differential equations with finite difference or finite element methods. Most of these methods deal with simple domains in the complex plane and cannot be applied on arbitrary regions. Furthermore, to the best of our knowledge, no work has been done on solving the Betrami equation on general Riemann surfaces. In this work, we are interested in developing an effective numerical algorithm to compute the quasi-conformal mapping on general Riemann surfaces of any genus. Of course, the developed algorithm could be easily applied to any arbitrary regions in the complex plane $\mathbb{C}$, since any regions in $\mathbb{C}$ are Riemann surfaces.

The fundamental problem in this paper is to find the quasi-conformal map $\phi$, associated to a given Beltrami differential. This can be done by solving the Beltrami equation using the proposed \emph{auxiliary metric method}. We firstly formulate a discrete analog of the quasi-conformal map on triangular meshes. Then, we propose an algorithm to compute the discrete quasi-conformal map. The basic idea is to construct a \emph{discrete auxiliary metric} based on the given Beltrami differential, such that $\phi$ becomes a conformal map under the auxiliary metric. We then use the Yamabe flow method to compute the conformal map under the new discrete metric. The resulting map is the desired quasi-conformal map. Numerically, the discrete quasi-conformal map converges to the continuous real solution as mesh size tends to 0.

The paper is laid out in the following way: Section II briefly review the most related works in the field; Section III introduces the theoretical background; Section IV describes the discrete Euclidean and hyperbolic Yamabe flow method; Section V explains how the Beltrami equations can be solved on general Riemann surfaces using the auxiliary metric. Section VI focuses on the computational methodologies; Section VII reports the experimental results; the paper is concluded in Section VIII.

\section{Previous work}
Conformal mapping has been broadly applied to surface parameterizations in digital geometry processing. Here we only review the most related works. We refer readers to the thorough surveys of \cite{floater2005,sheffer2007,HLS07} for various kinds of mesh parameterization techniques.

L{\'e}vy et al. in \cite{levy02} applied the Cauchy-Riemann equation for mesh parameterization and provided successful results on the
constrained 2D parameterizations with free boundaries. Desbrun et al. in \cite{desbrun2002} minimized the Dirichlet energy defined on
triangular meshes for computing conformal parameterization. Angle based flattening method (ABF) was introduced in
\cite{sheffer2001,sheffer2005}. Linearized version of ABF is introduced in \cite{ZLS07}. Gu and Yau in \cite{DBLP:conf/sgp/GuY03}
computed the conformal structure using Hodge theory. Gortler et al. in \cite{gortler2005} used discrete 1-forms for mesh parameterization with several holes. Ray et al. in \cite{ray2005} used the holomorphic 1-form to follow the principle curvatures for the quad remeshing purpose. K{\"a}lberer et al.
\cite{DBLP:journals/cgf/KalbererNP07} use branched covering to convert a given frame field on the surface to a vector field on the covering space. Spherical parameterizations are introduced in \cite{DBLP:journals/tog/GotsmanGS03, DBLP:journals/tog/PraunH03}. High genus surface parameterization is pioneered by Grimm and Hughes in \cite{cin03}. Recently, hyperbolic parameterization is introduced in \cite{JinTVCG08}.

Circle pattern was proposed by Bowers and Hurdal \cite{Bowers03}, and has been proven to be a minimizer of a convex energy by Bobenko and Springborn \cite{bobenko04}. An efficient circle pattern algorithm was developed by Kharevych et al. \cite{DBLP:journals/tog/KharevychSS06}. Discrete Ricci flow was introduced by Chow and Luo in \cite{ChowLuo03} and applied to graphics in \cite{JinTVCG08}. Ben-Chen et al. introduce an efficient method for scaling metrics to prescribed curvatures in \cite{DBLP:journals/cgf/Ben-ChenGB08}.

The theory for combinatorial Euclidean Yamabe flow was introduced by Luo in \cite{Luo04}. The theory for hyperbolic curvature flow was introduced
in \cite{BSP}. Springborn et al. \cite{Springborn08} identifies the Yamabe energy with the Milnor-Lobachevsky function and the heat
equation for the curvature evolution with the cotangent Laplace equation.

Recently, there has been various work on numerical quasi-conformal mapping techniques based on solving elliptic equations in the real plane with finite difference or finite element methods. Using finite difference methods to compute quasi-conformal maps on complex plane were proposed by Belinskii et al. \cite{Belinskii} and Mastin and Thompson \cite{Mastin}. These methods are difficult to implement for arbitrary regions. A finite difference scheme for constructing quasi-conformal mappings for arbitrary simply and doubly-connected region of the plane onto a rectangle was developed by Mastin and Thompson \cite{Mastin2}. Vlasynk \cite{Vlasyuk} applied similiar techniques for mappings of doubly connected and triply connected domains onto a parametric rectangle. A finite element based method was implemented by Weisel \cite{Weisel}. In \cite{Daripa2} Daripa proposed a numerical construction of quasi-conformal mappings in the plane using the Beltrami equation. The author presented an algorithm for the evaluation of one of the singular operators that arise in solving the Beltrami equation. The author subsequently applied the same method for numerical quasi-conformal mappings of exterior of simply connected domains onto the interior of a unit disk using the Beltrami equation \cite{Daripa}. This method was further extended to the quasi-conformal mapping of an arbitrary doubly connected domain with smooth boundaries onto an annulus $\Omega_R = \{\sigma: R< \sigma <1\}$ \cite{Daripa3}. All of these methods deal with simple domains in the complex plane and cannot be applied on arbitrary regions. Furthermore, as far as we know, no work has been done on solving the Betrami equation on general Riemann surfaces. In this work, we are interested in developing an effective numerical algorithm to compute the quasi-conformal mapping on general Riemann surfaces of any genus.

\section{Theoretical Background}
In this section, we briefly introduce the major concepts in differential geometry and Riemann surface theory, which are necessary to explain the quasi-conformal maps. We refer readers to \cite{diffgeombook,Kra04} for detailed information.

\subsection{Beltrami Equations and quasi-conformal maps}
Let $f: \mathbb{C}\to \mathbb{C}$ be a complex function. The following differential operators are more convenient for discussion
\[
    \frac{\partial}{\partial z} := \frac{1}{2}(\frac{\partial }{\partial x}-i\frac{\partial }{\partial
    y}),
    \frac{\partial}{\partial \bar{z}} := \frac{1}{2}(\frac{\partial }{\partial x}+i\frac{\partial }{\partial
    y}).
\]

$f$ is said to be {\it quasi-conformal} associated to $\mu$ if it is orientation-preserving and satisfies the following {\it Beltrami equation}:
\begin{equation}
 \frac{\partial f}{\partial \bar{z}} = \mu(z) \frac{\partial f}{\partial z}
\end{equation}

\noindent where $\mu(z)$ is some complex-valued Lebesgue measurable function satisfying $\sup|\mu| < 1$. $\mu$ is called the {\it Beltrami coefficient} of $f$. The Beltrami coefficient $\mu$ gives us all the information about the conformaity of $f$ (See Figure \ref{beltrami}).

If $\mu(z)= 0$ everywhere, $f$ is called \emph{holomorphic}. A holomorphic function satisfies the well-known Cauchy-Riemann equation
\[
    \frac{\partial f}{\partial \bar{z}}=0.
\]

Suppose $S$ is a surface embedded in $\mathbb{R}^3$, with the induced Euclidean metric $\mathbf{g}$. Let $U_\alpha \subset
S$ be an open set on $S$, with local parameterization is
$\phi_\alpha: U_\alpha \to \mathbb{C}$, such that the metric has local
representation
\[
    \mathbf{g} = e^{2\lambda(c)}dzd\bar{z},
\]
then $(U_\alpha,\phi_\alpha)$ is called an \emph{isothermal coordinate chart}. We can cover the whole surface by a collection of isothermal coordinate charts. All isothermal coordinate charts form a \emph{conformal structure} of the surface. The surface with a conformal structure is called a \emph{Riemann surface}.

Suppose $S_1$ and $S_2$ are two Riemann surfaces. $(U_\alpha,\phi_\alpha)$ is a local chart of $S_1$, $(V_\beta,\psi_\beta)$ is a local chart of $S_2$. $f:S_1\to S_2$ is a \emph{conformal map} if and only if
\[
  \psi_\beta\circ f\circ\phi_\alpha^{-1}:  \phi_\alpha(U_\alpha) \to \psi_\beta(V_\beta)
\]
is bi-holomorphic for all $\phi_\alpha$ and $\psi_\beta$. For simplicity, we still use $f$ to denote its local representation. Then a conformal map $f$ satisfies $\frac{\partial \phi}{\partial \bar{z}}=0$.

The definition of quasi-conformal maps of plane domains can be extended to Riemann surfaces. Instead of using the Beltrami coefficient, a global quantity called \emph{Beltrami differential} is used, which is independent of the choice of local parameters.

\begin{figure}
\centering
\includegraphics[height=2.5in]{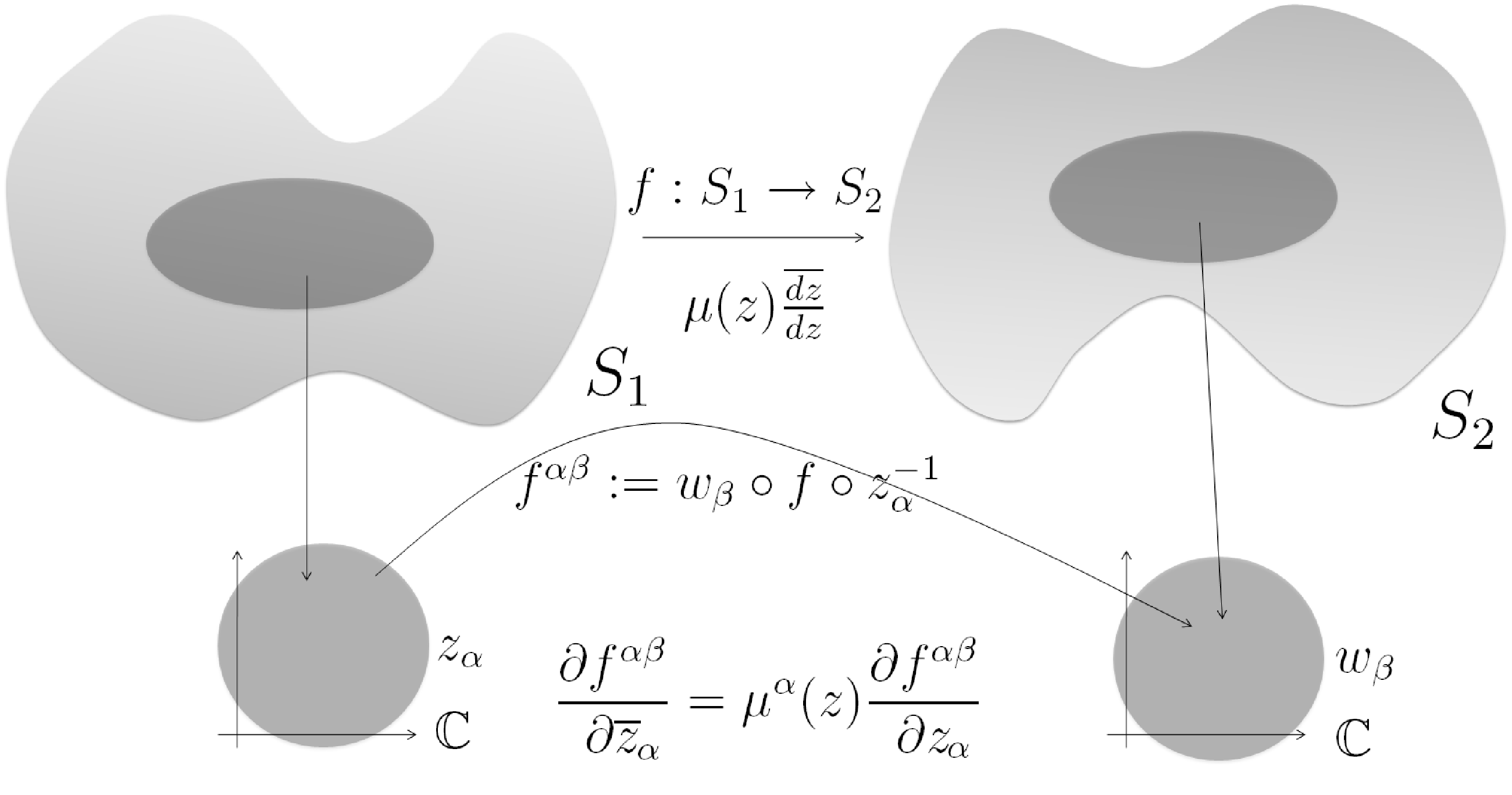}
\vspace{-2mm} \caption{The figure illustrates the definition of quasi-conformal maps between Riemann surfaces \label{fig:chart}}
\end{figure}

\bigskip

\begin{defn}[Beltrami Differentials]
A {\it Beltrami differential} $\mu(z)\frac{\overline{dz}}{dz}$ on a Riemann surface $R$ is an assignment to each chart $z_{\alpha}$ on $U_{\alpha}$ an $L_{\infty}$ complex-valued function $\mu^{\alpha}$ defined on $z_{\alpha}(U_{\alpha})$ such that:
\begin{equation} \label{beltramidifferential}
\mu^{\alpha}(z_{\alpha}) = \mu^{\beta}(z_{\beta}) \overline{(\frac{dz_{\beta}}{dz_{\alpha}})}/(\frac{dz_{\beta}}{dz_{\alpha}})
\end{equation}
\noindent on the region which is also covered by another chart $z_{\beta}$.
\end{defn}

\bigskip

Now, a quasi-conformal map between Riemann surfaces can be defined as follow (see Figure \ref{fig:chart} for the geometric illustration):

\bigskip

\begin{defn}[quasi-conformal maps between Riemann surfaces]
An orientation-preserving homeomorphism $f:S_1 \to S_2$ is called {\it quasi-conformal} associated to $\mu(z)\frac{\overline{dz}}{dz}$ if for any chart $z_{\alpha}$ on $S_1$ and any chart $w_{\beta}$ on $S_2$, the map $f^{\alpha \beta} := w_\beta\circ f\circ z_\alpha^{-1}$ is quasi-conformal associated with $\mu^{\alpha}(z_\alpha)$
\end{defn}

\bigskip

Note that the above definition is well defined. On a region of $S_1$ covered by two different charts $z_{\alpha}$ and $z_{\alpha '}$, we have
\begin{equation}
\mu^{\alpha '} = \frac{\partial f^{\alpha ' \beta}}{\partial \overline{z_\alpha'}}/\frac{\partial f^{\alpha ' \beta}}{\partial z_\alpha'} = (\frac{\partial f^{\alpha \beta}}{\partial \overline{z_\alpha}} \overline{\frac{dz_{\alpha}}{dz_{\alpha '}}})/(\frac{\partial f^{\alpha \beta}}{\partial z_\alpha} \frac{dz_{\alpha}}{dz_{\alpha '}}) = \mu^{\alpha}(z_{\alpha '}) \overline{(\frac{dz_{\alpha}}{dz_{\alpha '}})}/(\frac{dz_{\alpha}}{dz_{\alpha '}}) \nonumber
\end{equation}
This is guaranteed by Equation \ref{beltramidifferential}. Also, the definition does not depend on the chart $w_{\beta}$ used in the range of $f$. Let $w_{\beta}$ and $w_{\beta '}$ be two different charts on the range of $f$. We have
\begin{equation}
\begin{split}
\mu_\beta' (z_\alpha ) = \frac{\partial f^{\alpha \beta '}}{\partial \overline{z_\alpha}}/\frac{\partial f^{\alpha \beta '}}{\partial z_\alpha} &= (\frac{\partial w_{\beta '}}{\partial w_{\beta}}\frac{\partial f^{\alpha \beta}}{\partial \overline{z_{\alpha}}} + \frac{\partial w_{\beta '}}{\partial \overline{w_{\beta}}}\frac{\partial \overline{f^{\alpha \beta}}}{\partial \overline{z_{\alpha}}})/(\frac{\partial w_{\beta '}}{\partial w_{\beta}}\frac{\partial f^{\alpha \beta}}{\partial z_{\alpha}} + \frac{\partial w_{\beta '}}{\partial \overline{w_{\beta}}}\frac{\partial \overline{f^{\alpha \beta}}}{\partial z_{\alpha}})\nonumber \\
&= \frac{\partial f^{\alpha \beta }}{\partial \overline{z_\alpha}}/\frac{\partial f^{\alpha \beta }}{\partial z_\alpha} = \mu_\beta (z_\alpha )
\end{split}
\end{equation}
since $w_{\beta '}$ is holomorphic and so $ \frac{\partial w_{\beta '}}{\partial \overline{w_{\beta}}} = 0$.

\bigskip

With the above definitions, we can now formulate our problem of interest as follow:

\bigskip

\noindent $\mathbf{PROBLEM:}$ Let $S_1$ and $S_2$ be two Riemann surfaces with the same topology. Given a Beltrami differential $\mu(z)\frac{\overline{dz}}{dz}$, we are interested in finding the map $f$ such that for each chart $z_{\alpha}$ on $S_1$ and each chart $w_{\beta}$ on $S_2$, the map $f^{\alpha \beta} := w_\beta\circ f\circ z_\alpha^{-1}$ satisfies the partial differential equation:
\begin{equation}
\frac{\partial f^{\alpha \beta} }{\partial \overline{z_\alpha}} = \mu^{\alpha}(z_\alpha) \frac{\partial f^{\alpha \beta} }{\partial z_\alpha} \nonumber
\end{equation}

\subsection{Ricci Flow}
Let $S$ be a surface embedded in $\mathbb{R}^3$ with the induced
Euclidean metric $\mathbf{g}$. We say another Riemannian metric
$\bar{\mathbf{g}}$ is \emph{conformal} to $\mathbf{g}$, if there is
a scalar function $u:S\rightarrow\mathbb{R}$, such that
$\bar{\mathbf{g}}=e^{2u}\mathbf{g}$.

The Gaussian curvature induced by $\bar{\mathbf{g}}$ is
\[
\bar{K}=e^{-2u}(-\Delta_\mathbf{g}{u}+K), \label{eqn:curvature}
\]
where $\Delta_\mathbf{g}$ is the Laplace-Beltrami operator under
the original metric $\mathbf{g}$. The above equation is called the
\emph{Yamabe equation}. By solving the Yamabe equation, one can
design a conformal metric $e^{2u}\mathbf{g}$ by a prescribed
curvature $\bar{K}$.

Yamabe equation can be solved using \emph{Ricci flow} method. The
Ricci flow deforms the metric $\mathbf{g}(t)$ according to the
Gaussian curvature $K(t)$ (induced by itself), where $t$ is the time
parameter
\[
    \frac{dg_{ij}(t)}{dt} = 2(\bar{K}-K(t)) g_{ij}(t).
    \label{eqn:ricci_flow}
\]

The \emph{uniformization theorem} for surfaces says that any metric
surface admits a conformal metric, which induces constant Gaussian
curvature. The constant is one of $\{-1,0,+1\}$, determined by the
topology of the surface. Such metric is called the
\emph{uniformization metric}. Ricci flow converges to the
uniformization metric. Detailed proofs can be found in \cite{ric82}
and \cite{chow1991}.
\section{Discrete Euclidean and Hyperbolic Curvature Flow}
In practice, most surfaces are approximated by simplicial complexes, namely triangular meshes. Suppose $M$ is trianglar mesh, $V,E,F$ are vertex, edge and face set respectively. We use $v_i$ to denote the $ith$ vertex; $[v_i,v_j]$ the edge from $v_i$ to $v_j$; $[v_i,v_j,v_k]$ as the face, where the vertices are sorted counter-clock-wisely. On triangular meshes, we can derive a discrete version of Ricci flow, called the {\it discrete Yamabe flow}, which is analogous to the curvature flow on smooth surfaces. In this section, we describe in detail the discrete Euclidean and hyperbolic curvature flow that converge to the Euclidean and hyperbolic uniformization metric respectively.

On the discrete mesh, we can define the discrete metric, which is similar to the Riemannian metric. Basically, the discrete metric gives the length of each edge.

\medskip

\begin{defn}[Discrete Metric]
A discrete metric on a mesh $M$ is a function $l: E\to
\mathbb{R}^+$, such that on each triangle $[v_i,v_j,v_k]$, the
triangle inequality holds,
\begin{equation}
    l_{i} + l_{j} > l_{k}.\nonumber
\end{equation}
\end{defn}

\bigskip

\begin{figure}
\centering
\includegraphics[height=1.8in]{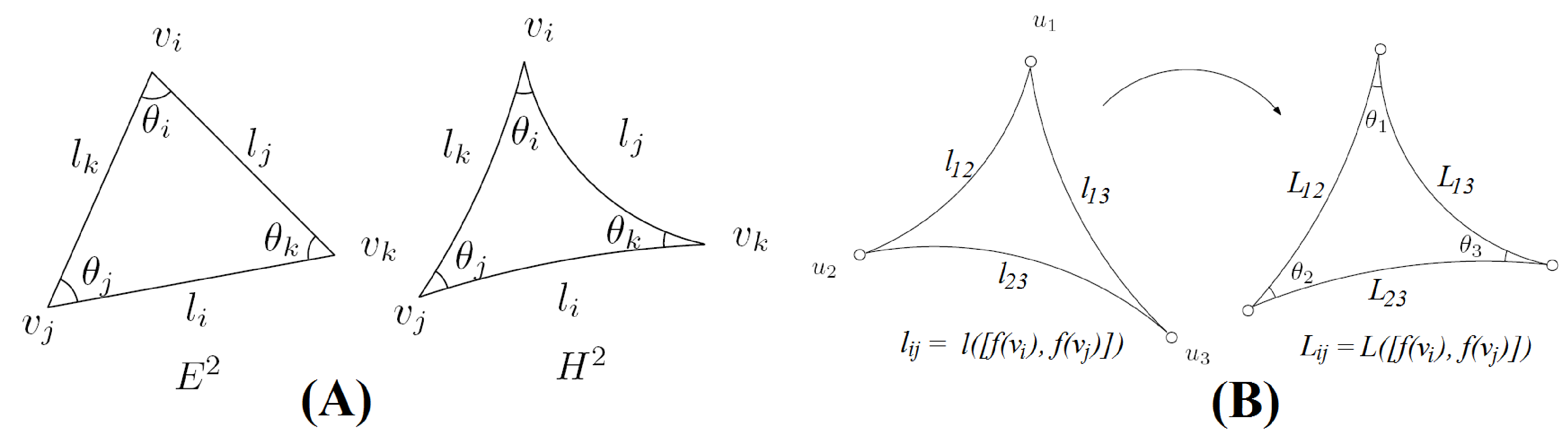}
\caption{(A) shows the Euclidean and Hyperbolic triangle. (B) shows the discrete surface Yamabe flow \label{fig:Yamabe}}
\end{figure}

Discrete metric represents a configuration of edge lengths. As shown
in figure \ref{fig:Yamabe}(A), different background geometries can be
assigned to a mesh.

\bigskip

\begin{defn}[Background Geometry]
Suppose $M$ is a mesh with a discrete metric. If all faces are Spherical (or Euclidean, or Hyperbolic triangle), then the mesh is with Spherical (or Euclidean, or Hyperbolic) background geometry, denoted as $\mathbb{S}^2$ (or $\mathbb{E}^2$, or $\mathbb{H}^2$).
\end{defn}

\bigskip

Discrete metric determines the corner angles on each face by the cosine law,
\begin{equation}
    \theta_i =
    \left\{
    \begin{array}{ll}
     \cos^{-1} \frac{l_j^2 + l_k^2 - l_i^2}{2l_jl_k} &
     \mathbb{E}^2\\
    \cos^{-1} \frac{\cosh l_j \cosh l_k - \cosh l_i}{2\sinh l_j \sinh l_k} &
    \mathbb{H}^2
    \end{array}
    \right.
\end{equation}

\bigskip

The discrete Gaussian curvature is defined as angle deficient,
\smallskip

\begin{defn}[Discrete Gaussian Curvature] Suppose $M$ is a
mesh with a discrete metric, which either in Euclidean or hyperbolic
background geometry. $[v_i,v_j,v_k]$ is a face in $M$,
$\theta_i^{jk}$ represent the corner angle at $v_i$ on the face. The
discrete Gaussian curvature of $v_i$ is defined as
\begin{equation}
    K_i = \left\{
    \begin{array}{rl}
    2\pi - \sum_{jk} \theta_i^{jk}& v_i \not\in \partial M\\
    \pi - \sum_{jk} \theta_i^{jk}& v_i \in \partial M\\
    \end{array}
    \right.
\end{equation}
\end{defn}

The total Gaussian curvature is controlled by the topology of the
mesh,
\begin{equation}
    \sum_i K_i + \lambda \sum_{ijk} A_{ijk} = 2\pi \chi(M),
\end{equation}
where $A_{ijk}$ is the area of the face $[v_i,v_j,v_k]$, $\lambda$
is zero if $M$ is with Euclidean background metric, $-1$ if $M$ is
with hyperbolic background metric.

\bigskip

\begin{defn}[Discrete conformal deformation]
Let $K$ be a triangulation mesh. Suppose $l$ and $L$ are two different discrete metrics on $K$. We say $L$ is a {\it discrete conformal deformation} of $l$ if there exists a function $u:V\to \mathbb{R}$, where $V$ is the set of all vertices of $K$, such that for all edges $[v_i,v_j]$ on $K_1$:

\begin{equation}\label{eqn:edgelength}
   L([f(v_i), f(v_j)]) =
    \left\{
    \begin{array}{ll}
    e^{u(v_i)} l([v_i, v_j]) e^{u(v_j)} & \mathbb{E}^2\\
    2 \sinh^{-1} e^{u(v_i} \sinh (\frac{l([v_i, v_j])}{2}) e^{u(v_j} &
    \mathbb{H}^2
    \end{array}
    \right.
\end{equation}
\noindent where $[v_i,v_j]$ is an edge on $K$. (\ref{fig:Yamabe}(B))
\end{defn}

\bigskip
Note that the definition of discrete hyperbolic conformality was due to Springborn et al. \cite{BSP}.  $u:V\to \mathbb{R}$ is called the discrete conformal factor. The discrete conformal factor is a function defined on every vertices of $K$.

\bigskip

\begin{defn}[Discrete Yamabe Flow]
Let $\overline{K}_i$ denote the user defined curvature at $v_i$, the
discrete Euclidean and Hyperbolic Yamabe flow has the same formula
\begin{equation}
\frac{du_i}{dt} = \overline{K}_i-K_i.
 \label{eqn:yamabe_flow}
\end{equation}
\end{defn}

\bigskip
The discrete Yamabe flow iteratively modifies the discrete conformal factor and hence define a sequence of conformal deformation. The initial metric is deformed into the uniformization metric.
\bigskip

In both Euclidean and hyperbolic cases, one can easily observe that under the discrete conformal deformation  (see equation \ref{eqn:thetaij}),
\begin{equation}
    \frac{\partial \theta_i}{\partial u_j}= \frac{\partial \theta_j}{\partial
    u_i}.
\end{equation}
Therefore, the differential 1-form $\sum_i K_i du_i$ is a closed 1-form. The following Yamabe energy is well defined.

\begin{figure*}[ht]
\centering
\includegraphics[height=1.8in]{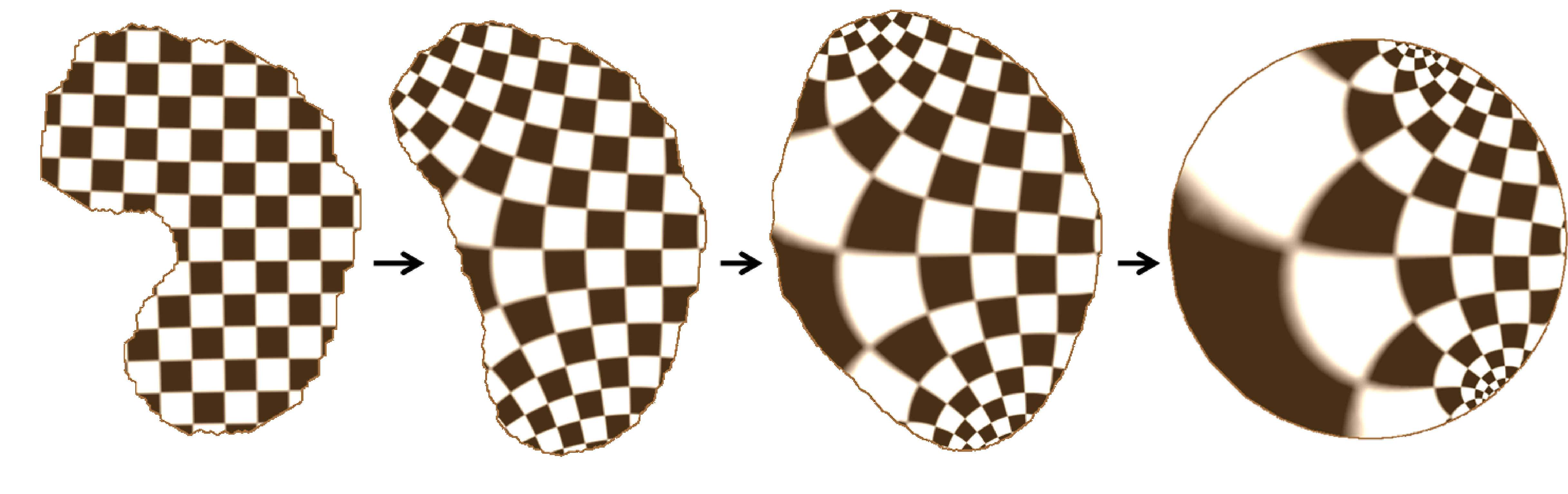}
\vspace{-2mm} \caption{The figure illustrates the idea of discrete curvature flow. A planar domain with non-zero curvature on the boundary is being deformed iteratively into a flat metric. Under the deformed metric, curvature is zero everywhere in the interior and is constant on the boundary (i.e. the boundary becomes a circle).  \label{fig:curvatureflowsc}}
\end{figure*}

\begin{defn}[Discrete Yamabe Energy]
Let $\mathbf{u}=(u_1,u_2,\cdots, u_n)$ be the conformal factor
vector, where $n$ is the number of vertices,
$\mathbf{u_0}=(0,0,\cdots, 0)$, then the discrete Euclidean and
Hyperbolic Yamabe energy has the same form as
\begin{equation}
E(\mathbf{u}) = \int_{\mathbf{u_0}}^{\mathbf{u}} \sum_{i=1}^n
(\overline{K}_i-K_i) du_i.
\end{equation}
\end{defn}

\bigskip

The discrete Yamabe flow is the negative gradient flow of the Yamabe
energy. The Hessian matrix $(h_{ij})$ of the discrete Euclidean and Hyperbolic Yamabe energy can be computed explicitly.

\bigskip
\begin{thm}
Let $[v_i,v_j]$ be an edge, connecting two faces $[v_i,v_j,v_k]$ and $[v_j,v_i,v_l]$, then Hessian matrix $(h_{ij})$ satisfies
\begin{equation} \label{eqn:hessian1}
    h_{ij} = \frac{\partial \theta_i^{jk}}{\partial u_j} + \frac{\partial \theta_i^{lj}}{\partial
    u_j}.
\end{equation}
and
\begin{equation}\label{eqn:hessian2}
    h_{ii} = \sum_{j,k} \frac{\partial \theta_i^{jk}}{\partial
    u_i},
\end{equation}
where the summation goes through all faces surrounding $v_i$, $[v_i,v_j,v_k]$.
Here,
\begin{equation}\label{eqn:thetaii}
    \frac{\partial \theta_i}{\partial u_i} =
    \left\{
    \begin{array}{ll}
    -\cot \theta_j - \cot \theta_k & \mathbb{E}^2\\
    -\frac{2c_ic_jc_k - c_j^2 - c_k^2 + c_i c_j + c_i c_k - c_j - c_k }{A(c_j+1)(c_k+1)} & \mathbb{H}^2\\
    \end{array}
    \right.
\end{equation}
\begin{equation} \label{eqn:thetaij}
    \frac{\partial \theta_i}{\partial u_j} =
    \left\{
    \begin{array}{lr}
    \cot \theta_k& \mathbb{E}^2\\
    \frac{c_i+c_j-c_k-1}{A(c_k + 1)} & \mathbb{H}^2\\
    \end{array}
    \right.
\end{equation}
where in $\mathbb{H}^2$ case, $c_k$ is $\cosh(y_k)$ and $A$ is
$\sin(\theta_k) \sinh(y_i)\sinh(y_j)$.
\end{thm}
\begin{proof}
See Appendix.
\end{proof}

\bigskip

By carefully examining the positive definiteness of
the Hessian, the convexity of the Yamabe energy can be obtained.

\bigskip

\begin{thm} The discrete Euclidean Yamabe energy is locally convex on the space of $\sum_i u_i = 0$. The discrete hyperbolic Yamabe energy is convex.
\end{thm}
\begin{proof}
See Appendix.
\end{proof}

\bigskip

The discrete Yamabe energy can be optimized using Newton's method directly. Given the mesh $M$, a conformal factor vector $\mathbf{u}$ is \emph{admissible}, if the deformed metric satisfies the triangle inequality on each face. The space of all admissible conformal factors is not convex in either Euclidean or Hyperbolic case. In practice, the step length in Newton's method needs to be adjusted. Once triangle inequality doesn't hold on a face, edge swap needs to be performed. After finitely many such surgery operations on the triangulation mesh, there will be no singularity developed in the normalized discrete Yamabe flow. We can then prove that the discrete Yamabe flow converges exponentially fast to the discrete metric with constant curvature. Specifically, we have

\bigskip

\begin{thm}
If no singularity develops in the discrete Yamabe flow after finitely many of surgery operations, the solution converges exponentially fast to a discrete metric with constant curvature as time approaches infinity. In other words,
\[
|K_i(t) - \overline{K}_i|\leq c_1 e^{-c_2 t}
\]
\noindent for some constants $c_1$ and $c_2$ and $K_i(t)$ is the discrete curvature at vertex $v_i$ at time $t$.
\end{thm}
\begin{proof}
See Appendix.
\end{proof}

\bigskip

Figure \ref{fig:convergence} shows the exponential convergence of the discrete Yamabe flow method. The Yamabe energy at each iteration of two different real human faces are plotted on the left and right respectively.

\begin{figure}
\centering
\includegraphics[height=2.2in]{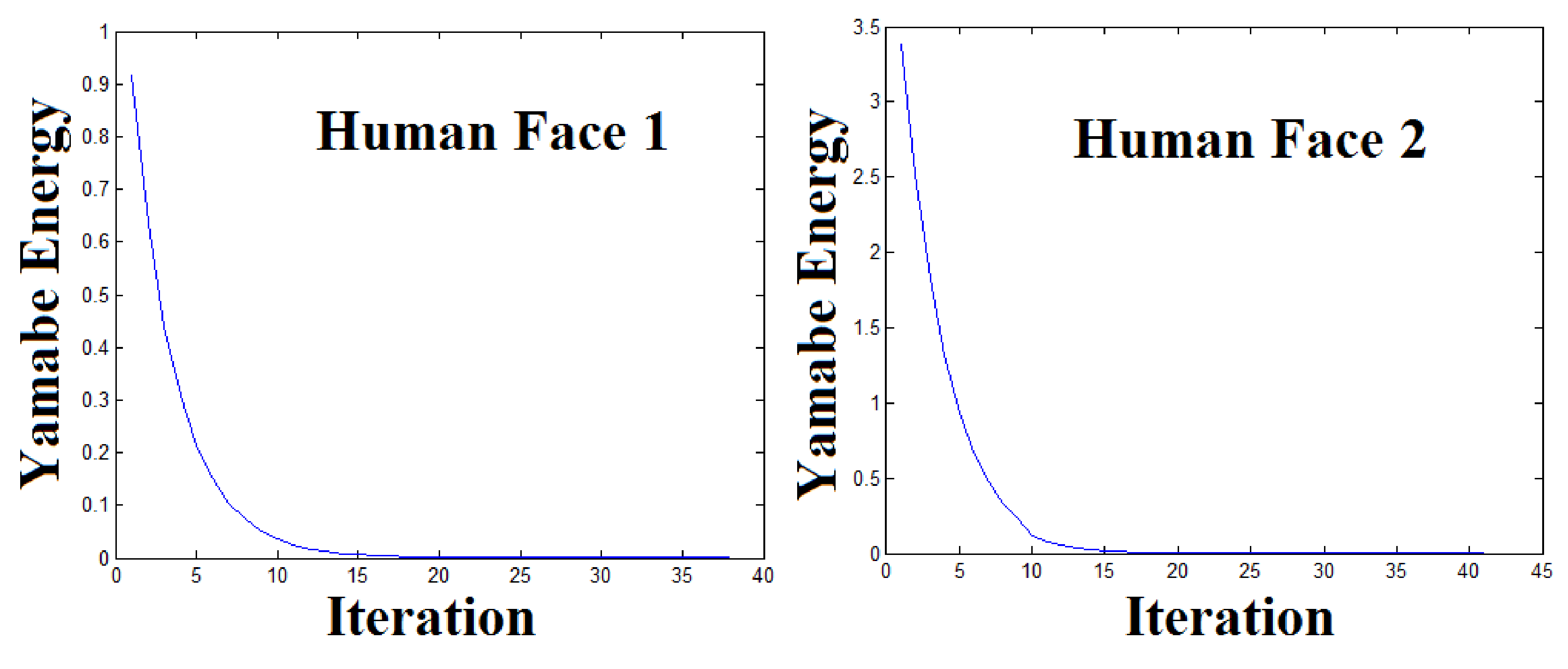}
\vspace{-2mm} \caption{The figure shows the exponential convergence of the discrete Yamabe flow method. The Yamabe energy at each iteration of two different real human faces are plotted on the left and right respectively.  \label{fig:convergence}}
\end{figure}

\section{Solving Beltrami Equations by Auxiliary Metric}
In this section, we prove the main theorem of this paper, which allows us to define an auxiliary metric to solve the Beltrami equation on the Riemann surface. It also lets us define the discrete version of quasi-conformal maps between meshes.

\bigskip

\noindent \begin{thm}[Auxiliary Metric associated to a Beltrami differential] \label{thm:auxiliary}
Suppose $(S_1,\mathbf{g}_1)$ and $(S_2,\mathbf{g}_2)$ are two metric surfaces, $\phi: S_1 \to S_2$ is a quasi-conformal map, the Beltrami differential is $\mu \frac{d\overline{z}}{dz}$. Let $z$ and $w$ be the local isothermal coordinates of $S_1$ and $S_2$ respectively, then $\mathbf{g}_1$ is $e^{2\lambda_1(z)}dzd\overline{z}$ and $\mathbf{g}_2$ is $e^{2\lambda_2(w)} dwd\overline{w}$. Define an auxiliary Riemannian metric on $S_1$,
\begin{equation}
    \mathbf{\tilde{g}_1} = e^{2\lambda_1(z)} |dz + \mu d\overline{z}|^2.
    \label{eqn:auxilary_metric}
\end{equation}
The auxiliary metric $\mathbf{\tilde{g}_1}$ is well-defined and the map $\phi:(S_1,\mathbf{\tilde{g}_1}) \to (S_2, \mathbf{g}_2)$ is a conformal map.
\end{thm}

\bigskip

\begin{proof}
We first prove the auxiliary metric $\mathbf{\tilde{g}_1}$ is well-defined. Consider the region which is covered by two different charts $z^{\alpha}$ and $z^{\beta}$. Suppose the local representations of $\mathbf{g}_1$ under $z^{\alpha}$ and $z^{\beta}$ are $e^{2\lambda_{\alpha}(z)}dz^{\alpha}d\overline{z^{\alpha}}$ and $e^{2\lambda_{\beta}(z)}dz^{\beta}d\overline{z^{\beta}}$ respectively.

Since $\frac{dz^{\alpha}}{d\overline{z^{\beta}}}= 0$, we have
\[
dz^{\alpha} = \frac{dz^{\alpha}}{dz^{\beta}}dz^{\beta} + \frac{dz^{\alpha}}{d\overline{z^{\beta}}}dz^{\overline{\beta}} = \frac{dz^{\alpha}}{dz^{\beta}}dz^{\beta}
\]
Also,
\[
e^{2\lambda_{\alpha}(z)}dz^{\alpha}d\overline{z^{\alpha}} = e^{2\lambda_{\alpha}(z)}|dz^{\alpha}|^2= e^{2\lambda_{\alpha}(z)}|\frac{dz^{\alpha}}{dz^{\beta}}|^2|dz^{\beta}|^2 = e^{2\lambda_{\beta}(z)}|dz^{\beta}|^2
\]
This gives, $e^{2\lambda_{\beta}(z)} = e^{2\lambda_{\alpha}(z)}|\frac{dz^{\alpha}}{dz^{\beta}}|^2$.

Thus,
\[
\begin{split}
e^{2\lambda_{\alpha}(z^{\alpha})}|dz^{\alpha}+ \mu^{\alpha} d\overline{z^{\alpha}}|^2 &=e^{2\lambda_{\alpha}(z^{\alpha})}|\frac{dz^{\alpha}}{dz^{\beta}}dz^{\beta}+ \mu^{\alpha} \overline{\frac{dz^{\alpha}}{dz^{\beta}}}d\overline{z^{\alpha}}|^2\\
&=e^{2\lambda_{\alpha}(z^{\alpha})}|\frac{dz^{\alpha}}{dz^{\beta}}|^2|dz^{\beta}+ \mu^{\alpha} (\overline{\frac{dz^{\alpha}}{dz^{\beta}}}/\frac{dz^{\alpha}}{dz^{\beta}})d\overline{z^{\alpha}}|^2\\
= e^{2\lambda_{\beta}(z^{\beta})}|dz^{\beta}+ \mu^{\beta} d\overline{z^{\beta}}|^2
\end{split}
\]
To see the map $\phi:(S_1,\mathbf{\tilde{g}_1}) \to (S_2, \mathbf{g}_2)$ is a conformal map, let $\phi^*\mathbf{g}_2$ denote the pull back metric,
\[
    \phi^*\mathbf{g_2} = e^{2\lambda_2( \phi(z) )} |d\phi(z)|^2
\]

Under the pull back metric, the map $\phi:( S_1, \phi^*\mathbf{g}_2) \to (S_2,\mathbf{g}_2)$ is isometric.
\[
\begin{array}{lcl}
    d\phi(z) &=& \frac{\partial \phi(z)}{\partial z} dz + \frac{\partial \phi(z)}{\partial \overline{z}} d \overline{z},\\
    &=& \frac{\partial \phi(z)}{\partial z} (dz + \mu d\overline{z} ).
\end{array}
\]
Therefore,
\[
    \phi^*\mathbf{g_2} = e^{2\lambda_2( \phi(z) )} |\frac{\partial \phi(z)}{\partial z}|^2  |dz+\mu d\overline{z}|^2 \label{eqn:eqn2}\\
\]
According to the definition of $\mathbf{\tilde{g}_1}$ in Equation
\ref{eqn:auxilary_metric}, $\phi^*\mathbf{g_2} = e^{2\lambda_2( \phi(z) ) - 2\lambda_1( z )} |\frac{\partial \phi(z)}{\partial z}|^2\mathbf{\tilde{g}_1}$. $\phi^*\mathbf{g_2}$ is conformal to
$\mathbf{\tilde{g}_1}$. Because $\phi:(S_1, \phi^*\mathbf{g_2})\to
(S_2,\mathbf{g}_2)$ is isometric, therefore $\phi:(S_1,
\mathbf{\tilde{g}_1})\to (S_2,\mathbf{g}_2)$ is conformal.
\end{proof}

\bigskip

This theorem tells us in order to solve the Beltrami equation on the Riemann surface, we simply need to define a new auxiliary metric associated with the prescribed Beltrami differential. We can then solve the Beltrami equation by computing a conformal map associated with the newly defined metric.
Specifically, we need to compute a conformal map $\phi:(S_1, \mathbf{\tilde{g}_1})\to (S_2,\mathbf{g}_2)$ (from Theorem 5.2). It can be done by flattening $(S_1, \mathbf{\tilde{g}_1})$ and $(S_2,\mathbf{g}_2)$ onto the unit disk. Let $\varphi_1: (S_1, \mathbf{\tilde{g}_1})\to \mathbb{D}$ and $\varphi_2: (S_2, \mathbf{\tilde{g}_2})\to \mathbb{D}$ be the conformal parameterizations of $S_1$ and $S_2$ respectively. The conformal map $\phi:(S_1, \mathbf{\tilde{g}_1})\to (S_2,\mathbf{g}_2)$ can be obtained by: $\phi = \varphi_2^{-1}\circ\varphi_1$. To compute the conformal parameterizations $\varphi_i$ ($i=1,2$), it suffices to find the flat metric of $(S_1, \mathbf{\tilde{g}_1})$ and $(S_2,\mathbf{g}_2)$. We can apply the discrete Yamabe flow method to conformally deform the original metric on the triangular mesh to the uniformization metric.

This motivates us to give a discrete definition of quasi-conformal maps between triangular meshes.

\bigskip

\begin{figure}
\centering
\includegraphics[height=2.2in]{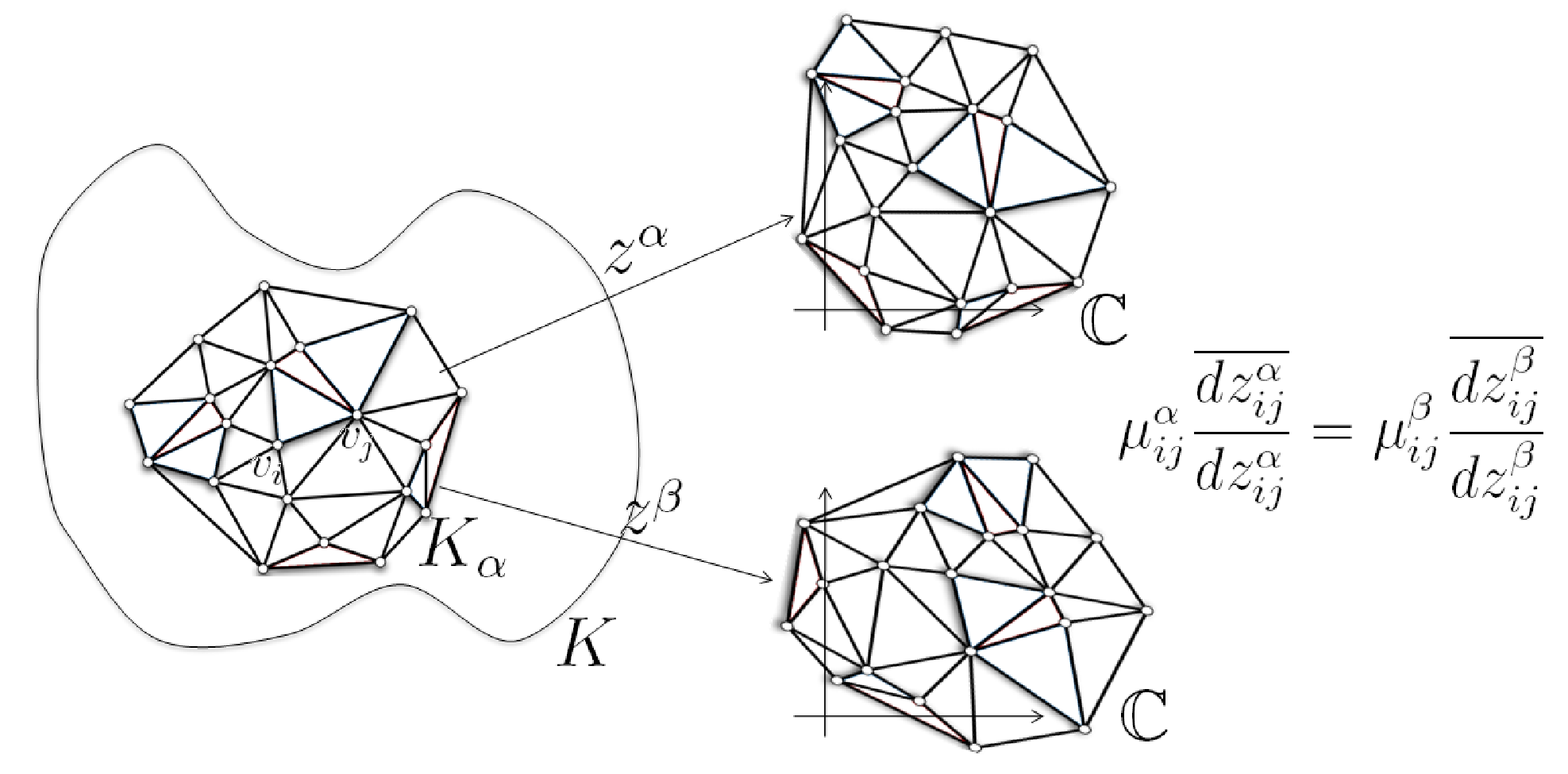}
\caption{ The figure illustrates the idea of discrete local chart and discrete Beltrami differential. \label{fig:localchart}}
\end{figure}

\begin{defn}[Discrete Local Charts]
Let $K$ be a triangular mesh. A mesh $K_{\alpha}$ is called a {\it submesh} of $K$ if every vertices, edges and faces of $K_{\alpha}$ belong to $K$. A {\it discrete local chart} $z^{\alpha}: K_{\alpha} \to z^{\alpha}(K_{\alpha}) \subset \mathbb{C}$ is a discrete conformal map from $K^{\alpha}$ to a mesh $z^{\alpha}(K^{\alpha})$ embedded in the complex plane (See Figure \ref{fig:localchart}).
\end{defn}

\bigskip

Since the triangular meshes we consider are discrete approximations of smooth surfaces, we can assume that the triangular meshes are covered by a collection of discrete local charts.

\bigskip

\begin{defn}[Discrete Beltrami Differential]\label{def:discretebeltrami}
A {\it discrete Beltrami differential} $\mu\frac{\overline{dz}}{dz}$ is an assignment to each discrete local chart $z^{\alpha}$ on $K_{\alpha}$ an $L_{\infty}$ complex-valued function $\mu^{\alpha}$ defined on $z^{\alpha}(K_{\alpha})$ such that:
\begin{equation}\label{eqn:discretebel}
(\frac{\mu^{\alpha}(v_i) + \mu^{\alpha}(v_j)}{2})\frac{\overline{z^{\alpha}(v_j)-z^{\alpha}(v_i)}}{z^{\alpha}(v_j)-z^{\alpha}(v_i)} = (\frac{\mu^{\beta}(v_i) + \mu^{\beta}(v_j)}{2})\frac{\overline{z^{\beta}(v_j)-z^{\beta}(v_i)}}{z^{\beta}(v_j)-z^{\beta}(v_i)}
\end{equation}
\noindent where $[v_i, v_j]$ is any edge in the region which is also covered by another chart $z_{\beta}$ (See Figure \ref{fig:localchart}).
\end{defn}

\bigskip

By letting $\mu^{\alpha}_{ij} = \frac{\mu^{\alpha}(v_i) + \mu^{\alpha}(v_j)}{2}$, $\mu^{\beta}_{ij} = \frac{\mu^{\beta}(v_i) + \mu^{\beta}(v_j)}{2}$, $dz^{\alpha}_{ij} = z^{\alpha}(v_j)-z^{\alpha}(v_i)$ and $dz^{\beta}_{ij} =z^{\beta}(v_j)-z^{\beta}(v_i)$, Equation \ref{eqn:discretebel} can be simplified as:

\[
\mu^{\alpha}_{ij}\frac{\overline{dz^{\alpha}_{ij}}}{dz^{\alpha}_{ij}} = \mu^{\beta}_{ij}\frac{\overline{dz^{\beta}_{ij}}}{dz^{\beta}_{ij}}
\]

\noindent which is an analog of Equation \ref{beltramidifferential}.

\bigskip

\begin{defn}[Discrete quasi-conformal Map]
Let $\mu\frac{d\overline{z}}{dz}$ be a given discrete Beltrami differential. A map $\phi: (K_1, l) \to (K_2, L)$ between meshes $K_1$ and $K_2$ is called {\it discrete quasi-conformal} if: with respect to a new metric $\widetilde{l}$ on $K_1$, the map $\phi:(K_1, \widetilde{l})\to (K_2, L)$ is discrete conformal where:

\noindent $\widetilde{l}([v_i, v_j]) := l([v_i, v_j])\frac{|dz_{ij} + \mu_{ij} \overline{dz_{ij}}|}{|dz_{ij}|}$; $dz_{ij} = z(v_j) - z(v_i)$; $\mu_{ij} = \frac{\mu_i + \mu_j}{2}$ for any local isothermal coordinates $z$ of $K_1$.
\end{defn}

\bigskip
$\widetilde{l}$ is called the {\it discrete auxiliary metric} associated with $\mu\frac{d\overline{z}}{dz}$. Note that the definition is well-defined. Suppose an edge $[v_i, v_j]$ is covered by both charts $z^{\alpha}$ and $z^{\beta}$, we have
\[
\begin{split}
l([v_i, v_j])\frac{|dz^{\alpha}_{ij} + \mu^{\alpha}_{ij} \overline{dz^{\alpha}_{ij}}|}{|dz^{\alpha}_{ij}|} &= l([v_i, v_j])|1+\mu^{\alpha}_{ij}\frac{\overline{dz^{\alpha}_{ij}}}{dz^{\alpha}_{ij}}| \\
& = l([v_i, v_j])|1+\mu^{\beta}_{ij}\frac{\overline{dz^{\beta}_{ij}}}{dz^{\beta}_{ij}}| \\
& = l([v_i, v_j])\frac{|dz^{\beta}_{ij} + \mu^{\beta}_{ij} \overline{dz^{\beta}_{ij}}|}{|dz^{\beta}_{ij}|}
\end{split}
\]
\bigskip

After the discrete auxiliary metric $\widetilde{l}$ is computed for each edge, we use the discrete Yamabe flow method introduced in Section IV to deform the metric into the uniformization metric and obtain the quasi-conformal map associated with the given Beltrami differential.

\section{Computational Algorithms}
In this section, we describe in detail the algorithms for Euclidean and Hyperbolic Yamabe flow. We also describe the algorithm for computing the discrete auxiliary metric associated uniquely with the prescribed discrete Beltrami differential.

\bigskip

\subsection{Yamabe Flow}
Algorithm 1 describes the details for Euclidean and Hyperbolic Yamabe flow.

\bigskip

\begin{center}
\underline{$\mathbf{Algorithm\ 1:}$ Euclidean $\&$ Hyperbolic Yamabe Flow}
\end{center}

\smallskip

\noindent  $\mathbf{Require:}$ A triangular mesh $M$, the target curvature $\overline{K}$.

\begin{enumerate}
\item Initialize $u_i$ to zero for each vertex $v_i$.
\item $\mathbf{Repeat}$
\item For each edge, compute the edge length using formula \ref{eqn:edgelength}.
\item For each face $[v_i,v_j,v_k]$ compute the corner angles $\theta_i,\theta_j$ and $\theta_k$.
\item For each face $[v_i,v_j,v_k]$, compute $\frac{\partial\theta_i}{\partial u_j}, \frac{\partial \theta_j}{\partial u_k}$, and $\frac{\partial\theta_k}{\partial u_i}$ using formula \ref{eqn:thetaij}.
\item For each face $[v_i,v_j,v_k]$, compute $\frac{\partial\theta_i}{\partial u_i},\frac{\partial \theta_j}{\partial u_j}$ and $\frac{\partial\theta_k}{\partial u_k}$ using formula \ref{eqn:thetaii}.
\item Construct the Hessian matrix  $H$ using formula \ref{eqn:hessian1} and \ref{eqn:hessian2}.
\item Solve linear system $H\delta \mathbf{u} = \overline{K}-K$.
\item Update discrete conformal factor $\mathbf{u} \leftarrow \mathbf{u}+\delta \mathbf{u}$.
\item For each vertex $v_i$, compute the Gaussian curvature $K_i$.
\item $\mathbf{until}$ $\max_{v_i\in M}|\overline{K}_i - K_i| < \epsilon$
\end{enumerate}

\subsection{Embedding}

After computing the discrete metric of the mesh, we can embed the
mesh onto $\mathbb{R}^2$ or $\mathbb{H}^2$. Here, we discuss the
hyperbolic case. Euclidean case is very similar.

In this work, we use Poincar\'{e} disk to model the hyperbolic space
$\mathbb{H}^2$, which is the unit disk $|z|<1$ on the complex plane
with the metric
\[
    ds^2 = \frac{4dzd\overline{z}}{(1-z\overline{z})^2}.
\]
The rigid motion is the M\"obius transformation
\[
    z \to e^{i\theta} \frac{z-z_0}{1-\overline{z}_0z},
\]
where $\theta$ and $z_0$ are parameters. A hyperbolic circle
$(\mathbf{c},r)$ is also a Euclidean circle $(\mathbf{C},R)$ with
\[
    \mathbf{C}=\frac{2-2\mu^2}{1-\mu^2\mathbf{c}\overline{\mathbf{c}}}
    \mathbf{c};
\]
\[
 R^2 = \mathbf{C}\overline{\mathbf{C}} -
    \frac{\mathbf{c}\overline{\mathbf{c}}-\mu^2}{1-\mu^2\mathbf{c}\overline{\mathbf{c}}},
\]
where $\mu = \tanh \frac{r}{2}$. Given a two points $p$ and $q$ on $\mathbb{H}^2$, the unique geodesic through them is a circular arc joining them and is perpendicular to the unit circle.

Basically, we can isometrically flatten triangle by triangle using the hyperbolic cosine
law, as described in algorithm 2. .

\bigskip

\begin{center}
\underline{$\mathbf{Algorithm\ 2:}$ Embed on $\mathbb{H}^2$ }
\end{center}

\smallskip

\noindent  $\mathbf{Require:}$  A triangular mesh $M$, a set of fundamental group generators intersecting only at the base point $p$, using the algorithm in \cite{JinTVCG08}

\begin{enumerate}
\item Slice $M$ along the base loops to form a fundamental domain
$\overline{M}$.
\item Embed the first triangle $[v_0,v_1,v_2]\in\overline{M}$,
\[
    \tau(v_0)=(0,0), \tau(v_1) = \tanh \frac{l_{01}}{2}, \tau(v_2) =
    \tanh \frac{l_{02}}{2}e^{i\theta_0^{12}}.
\]
\item Put all the neighboring faces of the first face to a face
queue.

\medskip

$\mathbf{While:}$ the face queue is not empty
\item Pop the first face from
the queue $[v_i,v_j,v_k]$

\item Suppose $v_i$ and $v_j$ has been embedded, compute the
intersection of two hyperbolic circles
\[
(\tau(v_i),l_{ik}) \cap (\tau(v_j),l_{jk})
\]
by converting them to Euclidean circles,
\item $\tau(v_k)$ is
chosen the keep the orientation of the face upward.
\item Put the neighboring faces, which haven't accessed, to the queue.
\end{enumerate}

\subsection{Auxiliary Metric}

The following algorithm 3 computes the discrete auxiliary metric based on the discrete Beltrami differential.

\bigskip

\begin{center}
\underline{$\mathbf{Algorithm\ 3:}$ Auxiliary Metric}
\end{center}

\smallskip

\noindent  $\mathbf{Require:}$ A triangular mesh $M$, a conformal parameterization $z: V\to \mathbb{C}$, discrete Beltrami differential $\mu :V \to \mathbb{C}$.

\begin{enumerate}
\item For each edge $[v_i,v_j]$, compute the edge length $l_{ij}$ using the Euclidean metric \\
$\mathbf{For\ all:}$ Edge $[v_i,v_j]$ $\mathbf{do}$
\item $dz \leftarrow z(v_j) - z(v_i)$
\item $\mu \leftarrow \frac{1}{2} ( \mu( v_i ) + \mu( v_j) )$
\item $\lambda \leftarrow \frac{|dz + \mu d\overline{z}|}{|dz|}$
\item $l_{ij} \leftarrow \lambda l_{ij}$ \\
$\mathbf{end\ for}$
\end{enumerate}

\section{Experimental results}
We implement our algorithm using generic C++
on Windows platform. The linear systems are solved using conjugate
gradient method. The experiments are carried out on a laptop with
2.0 GHZ CPU, 3.00 G RAM. The human face surfaces are captured using
phase shifting structured light method. Computational time is
reported in Table \ref{tab:time}.

\subsection{Quasi-conformal Map for Genus Zero Surfaces}

\begin{figure}
\centering
\includegraphics[height=5in]{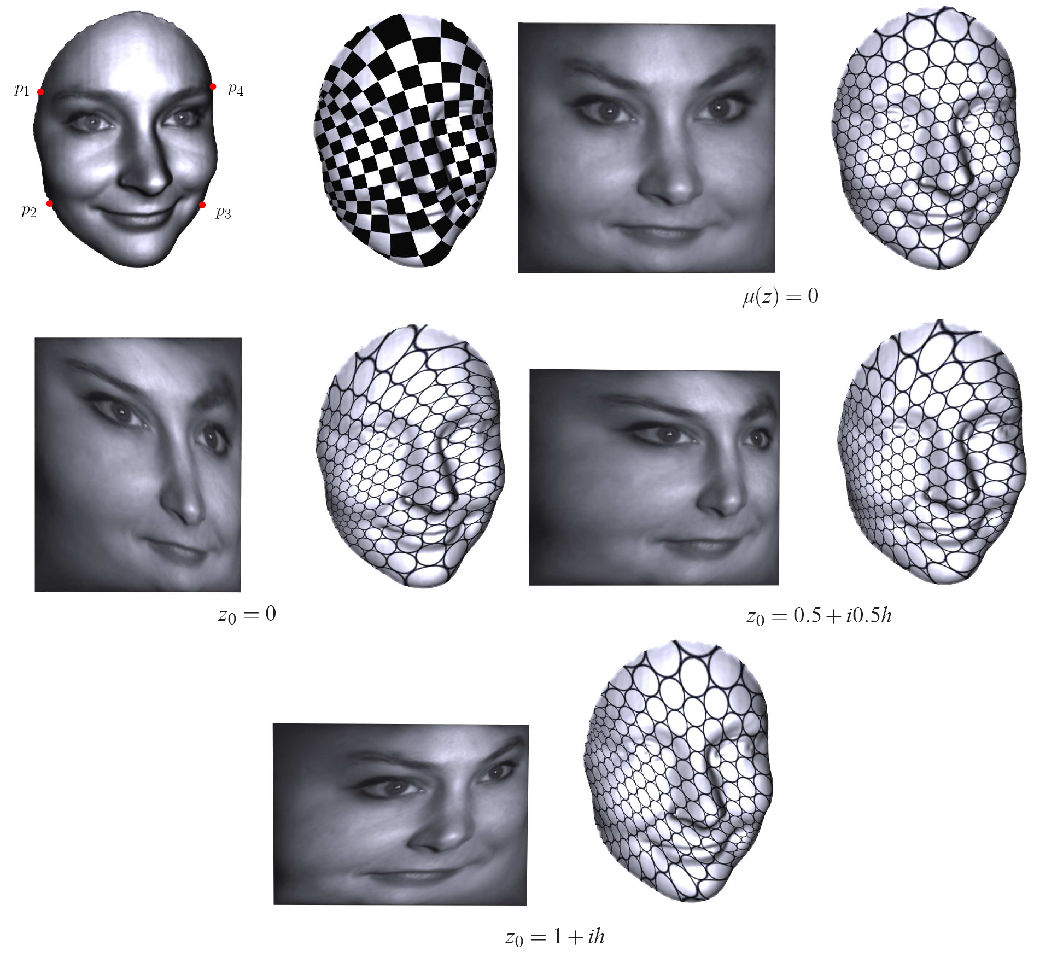}
\caption{Conformal and Quasi-conformal maps for a face
surface. The conformal parameter domain is a rectangle
with unit width and $h$ height on the 1st row. For all other
rows $\mu = \frac{z-z_0}{2\sqrt{1+h^2}}$with different $z_0$'s.  \label{fig:face4points}}
\end{figure}

Figure \ref{fig:face4points} shows the experimental results for quasi-conformal maps of the human face surface. The original face is shown in the top left corner. Four corner vertices are selected on the boundary, shown as $p_0,p_1,p_2,p_3$. We set the target curvature to be $\frac{\pi}{2}$ for those corner vertices, and zero for all vertices everywhere else. The Yamabe flow conformally maps the surface onto a planar rectangle. The corner vertices are mapped to the rectangle corners. We set the left lower corner to be the origin, the edges to be parallel to the axes, the
width to be $1$. Then the height $h$ gives us the conformal module of the original face surface with four fixed corners (topological quadrilateral). This provides us the conformal parameter $z$ of the surface.

In the second and the third row, we set different Beltrami coefficients. The image of the quasi-conformal map is shown on the left, the circle packing texture mapping is shown on the right. The Beltrami coefficient is set to be $\mu = \frac{z-z_0}{2\sqrt{1+h^2}}$, with different values of $z_0$ for different cases. It is obvious that, the conformal module of the surface changes with different Beltrami coefficients.

\begin{figure}
\centering
\includegraphics[height=4in]{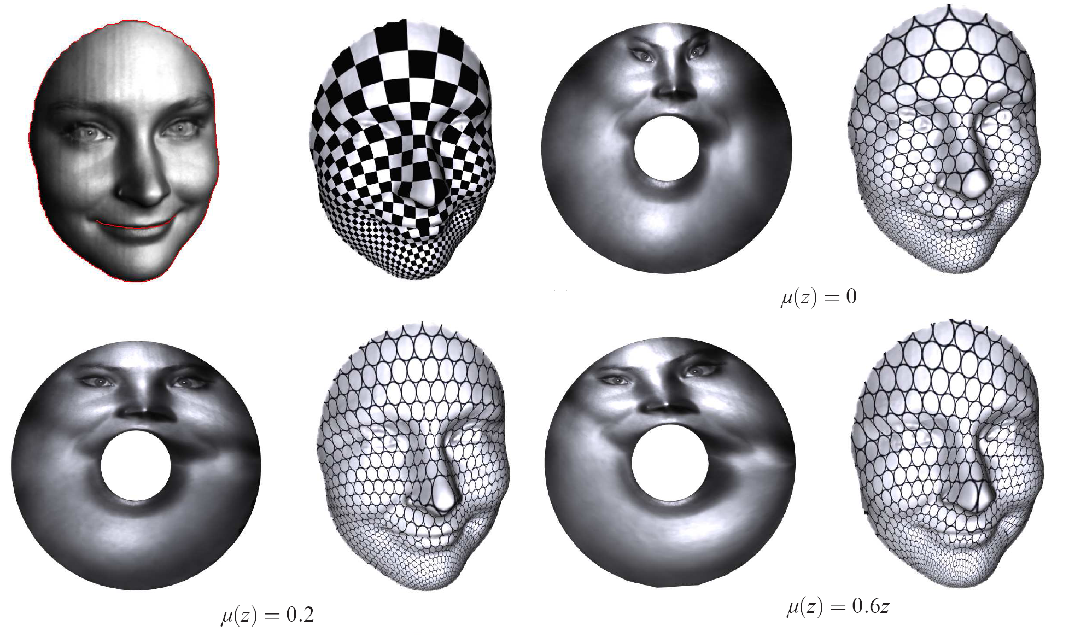}
\caption{Conformal and Quasi-conformal maps for a topological annulus. \label{fig:facemultiply1}}
\end{figure}

\begin{figure}
\centering
\includegraphics[height=5in]{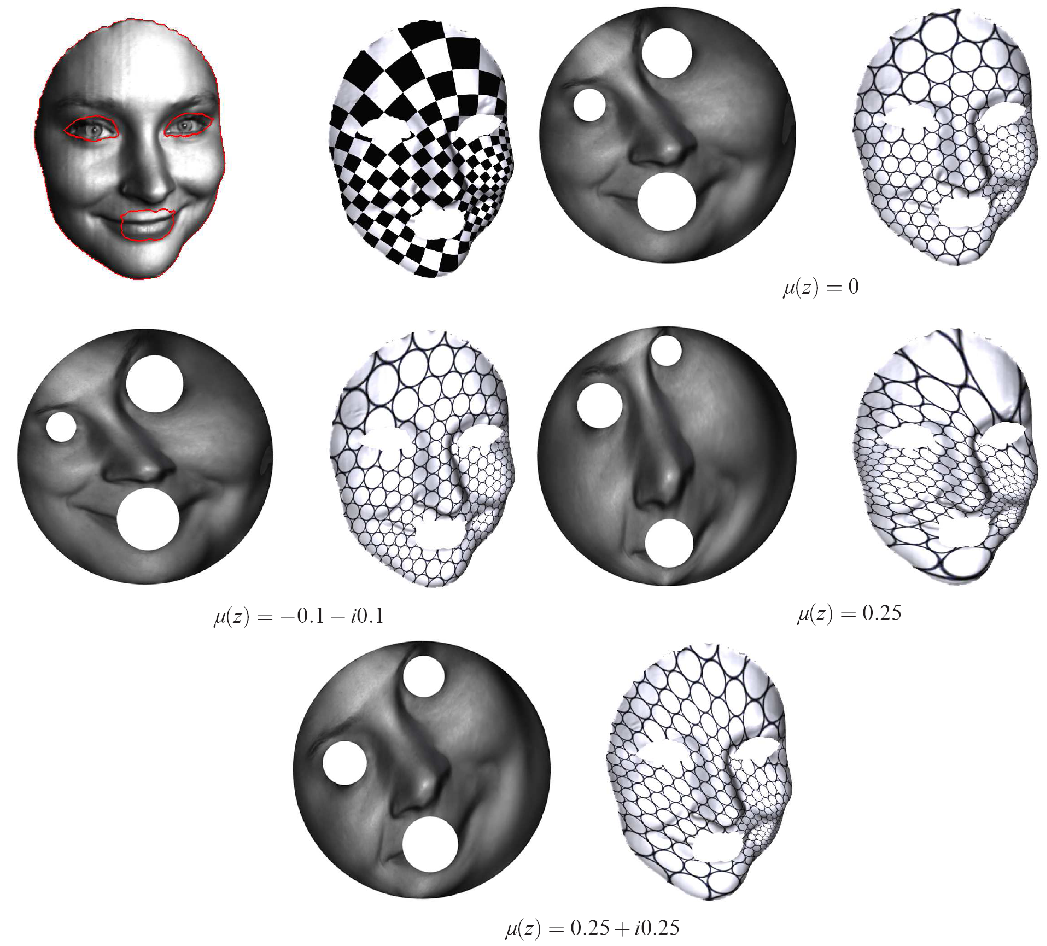}
\caption{Conformal and Quasi-conformal maps for a multiply connected domain. \label{fig:facemultiply3}}
\end{figure}

Figure \ref{fig:facemultiply1} and \ref{fig:facemultiply3} show the quasi-conformal maps for genus zero surfaces with multiply holes. In Figure \ref{fig:facemultiply1}, the human face surface is sliced open along the lip of the mouth which results in a doubly-connected open surface. Again, we set different Beltrami coefficients and compute the associated discrete auxiliary metrics. Using the discrete Yamabe flow, we conformally map the surface onto the annulus with respect to different auxiliary metrics. The target curvature is set to be zero in the interior and constant along the boundaries. The radius of the inner circles are different with different Beltrami coefficients, indicating a change in the conformal module. Figure \ref{fig:facemultiply3} shows the similar results for the genus zero human face surface with three slices (topological disk with 3 holes).

\begin{figure}
\centering
\includegraphics[height=5in]{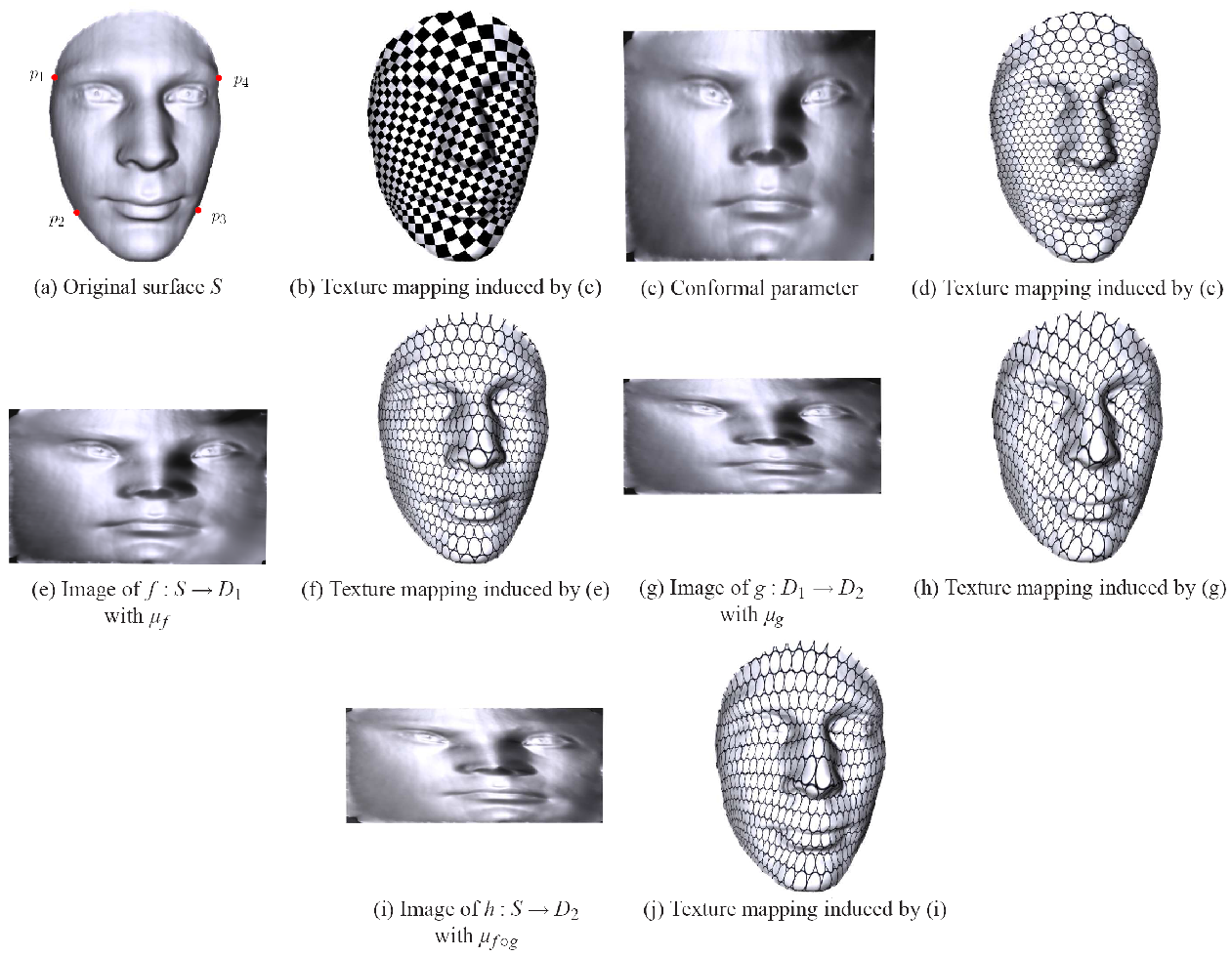}
\caption{Composed quasi-conformal maps for a topological quadrilateral.  \label{fig:face4pointscompose}}
\end{figure}

\begin{figure}
\centering
\includegraphics[height=6in]{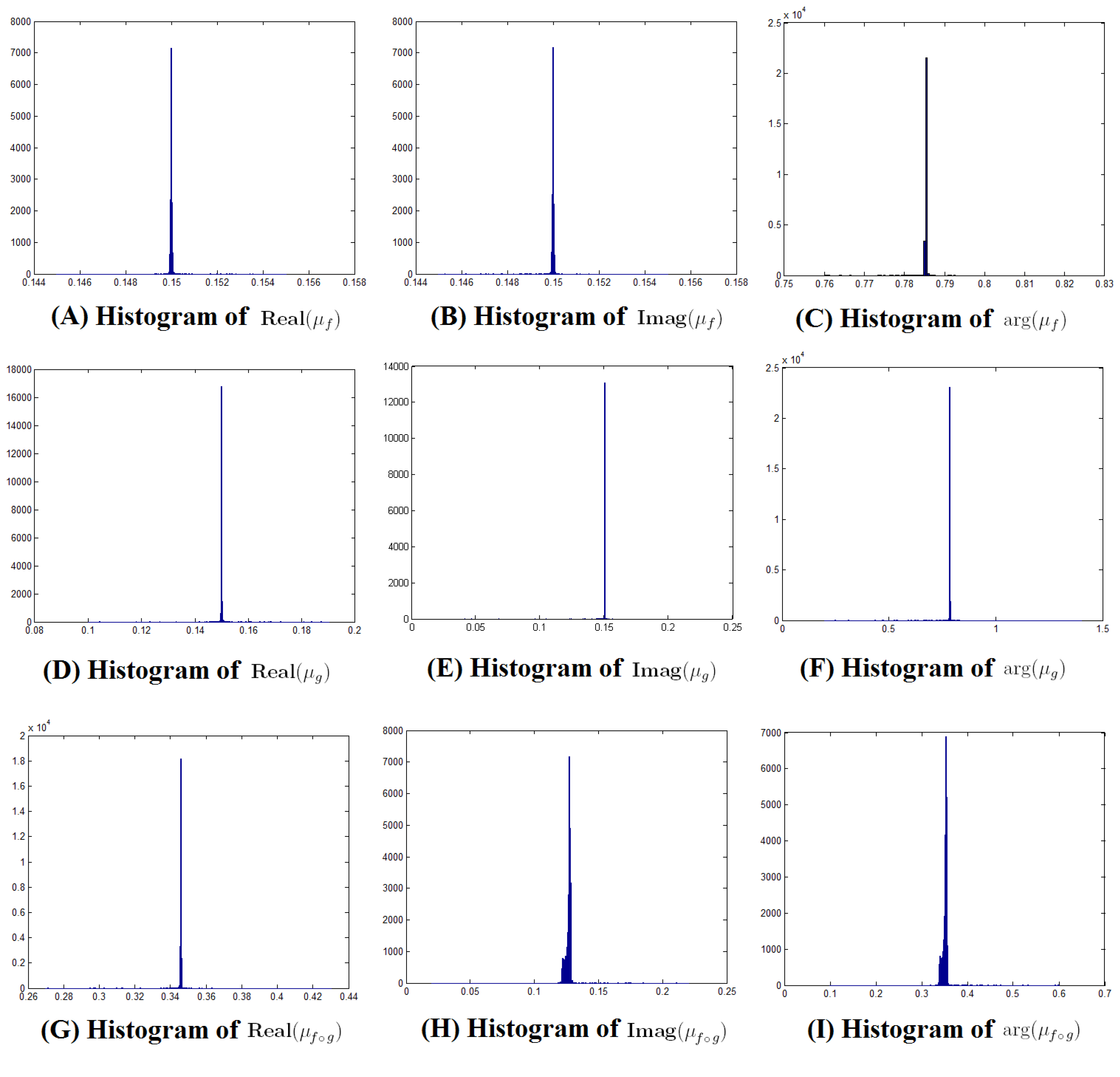}
\caption{Histogram of the real part, imaginary part and argument of the Beltrami coefficients.  \label{fig:histogram}}
\end{figure}

\subsection{Composition of Quasi-conformal Maps}
In the following experiment, we test the accuracy of our algorithm by computing the composed quasi-conformal maps using different approaches, and comparing their difference. If our method is accurate, the difference between the results obtained from the two approaches should be small.

Let $f: S\to D_1$ be a quasi-conformal map with Beltrami coefficient $\mu_f$, $g: D_1 \to D_2$ with Beltrami coefficient $\mu_g$. Then
the composed map $g\circ f: S \to D_2$ should have the Beltrami coefficient
\begin{equation}
    \mu_{g\circ f} = \frac{\mu_f + (\mu_g \circ f) \tau }{1 + \overline{\mu_f} (\mu_g\circ f)
    \tau},
    \label{eqn:mu_compose}
\end{equation}
where $\tau = \frac{\overline{f_z}}{f_z}$.

As shown in figure \ref{fig:face4pointscompose}, in our experiment, the original surface is a human face surface with four corner points (a
topological quadrilateral), as shown in (a). We compute its conformal parameter domain, as shown in (b), (c) and (d). Then we set the Beltrami coefficient $\mu_f$ to be $0.15+i0.15$, and use our method to compute a quasi-conformal map $f: S\to D_1$. The mapping result of $f$ is shown in (e) and (f). We set $\mu_g=0.15+i0.15$, and compute the quasi-conformal map $g:D_1 \to D_2$, as shown in (g) and (h).

We use the formula in Equation \ref{eqn:mu_compose} to compute the Beltrami coefficient for the composed map $\mu_{g\circ f} =0.34+i0.12$. We then solve the Beltrami equation $h_{\overline{z}} = \mu_{g\circ f} h_z$ to get a quasi-conformal map $h: S\to D_2$. In theory, $h$ should coincide with $g\circ f$. Our experimental result shows that $h$ is consistent with $g\circ f$. By comparing the result in (g) and (i), we can see the results of $g\circ f$ and $h$ are almost identical. We further measure the deviation between them numerically, using the following formula,
\[
    d(f,g) = \frac{1}{diag(S)A}\int_S |f(p)-g(p)| dp,
\]
where $A$ is the area of $S$, $diag(S)$ is the diagonal of the bounding box of $S$. The distance is the $L^1$ norm between $f$ and $g$, normalized by the diagonal of surface. In our experiment, the distance is $0.000044$, which is very small. This shows that our quasi-conformal map method is accurate. Figure \ref{fig:histogram} shows the histogram of the real part, imaginary part and argument of $\mu_f$, $\mu_g$ and $\mu_{f\circ g}$. (A), (B) and (C) shows the histograms of the real part, imaginary part and argument of the Beltrami coefficient $\mu_f$ of $f$ computed by our method. The histograms show that $\mathbf{Real}(\mu_f)=0.15$, $\mathbf{Imag}(\mu_f)=0.15$ and $\arg(\mu_f)=0.7854$ on almost all vertices, which agree with the exact solution. (D), (E) and (F) shows the histograms of the Beltrami coefficient $\mu_g$ of $g$. (G), (H) and (I) shows the histogram of the real part, imaginary part and argument of the Beltrami coefficient $\mu_{f\circ g}$ of the composition map $f\circ g$. The histograms show that $\mathbf{Real}(\mu_{f\circ g})=0.35$, $\mathbf{Imag}(\mu_f)=0.12$ and $\arg(\mu_f)=0.34$ on almost all vertices, which agree with the exact solution.

\subsection{Quasi-conformal Maps for Genus One Closed Surface}
\begin{figure*}[ht]
\centering
\includegraphics[height=5.2in]{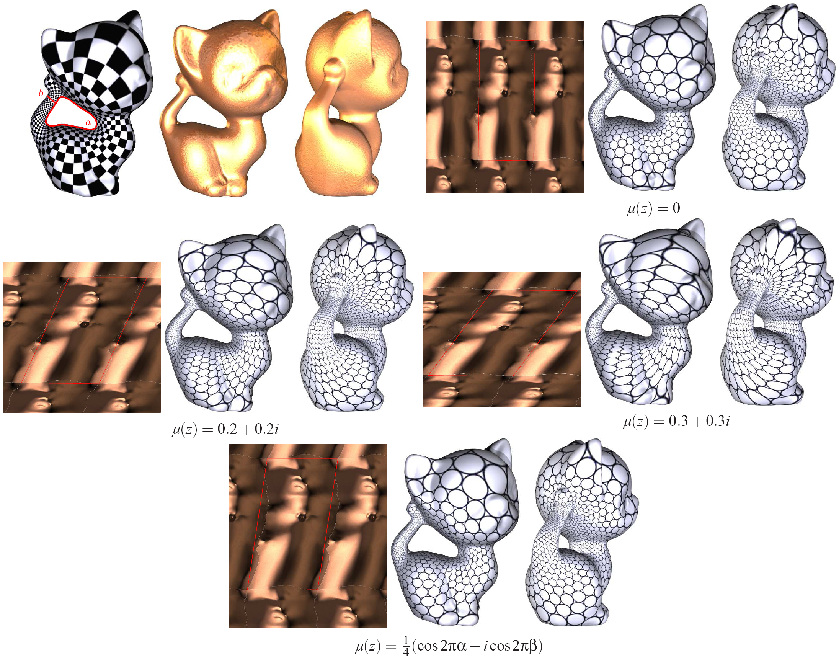}
\caption{Conformal and Quasi-conformal maps for a genus one surface.\label{fig:kitty}}
\end{figure*}

We test our algorithm for genus one closed surface as shown in Figure \ref{fig:kitty}. We first set $\mu$ to be zero, and compute a conformal flat metric using the curvature flow. Then we compute a homology basis, $\{a,b\}$ as shown in the leftmost frame of the top row. We then embed a finite portion of the universal covering space of the Kitten surface using the flat metric, shown in the right frame of the first row. The red rectangle shows the fundamental polygon, which is a parallelogram, with two adjacent edges $z_a$ and $z_b$.
The lattice $\Gamma$ is formed by the translations generated by
$z_a$ and $z_b$,
\[
    \Gamma = \{m z_a + n z_b | m,n \in \mathbb{Z}\}.
\]
The Kitten surface can be represented as the quotient space $M=\frac{\mathbb{R}^2}{\Gamma}$. This gives the conformal parameter domain of the surface. The rightmost frame of the first row illustrates the circle packing texture mapping induced by the conformal parameterization. In the second and the third row, we set different Beltrami coefficients. $\mu(z)$ are constants for the second row.

For the last row, the Beltrami coefficient is defined in a more complicated way. Because $\mu$ is defined on the Kitten surface, then it must satisfy the following consistency condition $\mu(z) = \mu( z + m z_a + n z_b )$. Given a point $z\in \mathbb{C}$, we can find a pair of real numbers $\alpha,\beta \in [0,1)$, such that $z\equiv \alpha z_a + \beta z_b ( mod \Gamma )$. Then $\mu$ is defined as $\mu(z) = \frac{1}{4}(\cos 2\pi \alpha + i \cos 2\pi \beta )$, which satisfies the above consistency condition.

\subsection{Quasi-conformal Maps for Genus Two Surfaces}
\begin{figure*}[ht]
\centering
\includegraphics[height=5.5in]{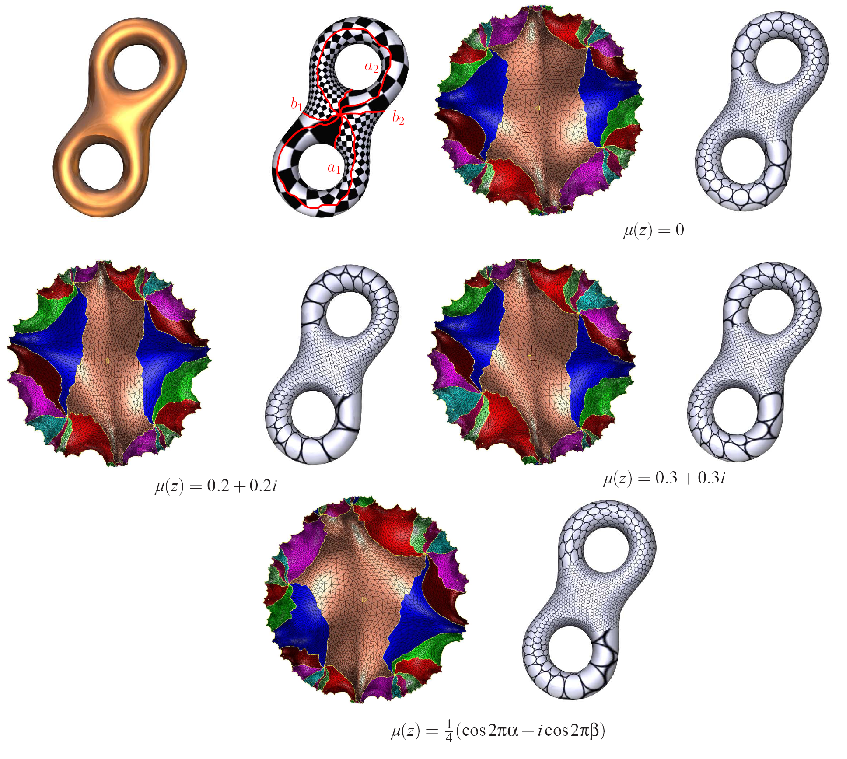}
\caption{Conformal and Quasi-conformal maps for a genus two surface.\label{fig:torus2}}
\end{figure*}

Our method can compute quasi-conformal maps for high genus surfaces, as shown in Figure \ref{fig:torus2}. We use hyperbolic Yamabe flow to compute the hyperbolic metric of the surface, then the homology basis $\{a_1,b_1,a_2,b_2\}$. We flatten a fundamental domain, compute the Fuchsian group generators, and flatten a finite portion of the universal covering space of the surface. Details can be found in \cite{cin03,JinTVCG08}. This gives a conformal atlas of the surface. Because of the difference between
hyperbolic metric and Euclidean metric, the texture mappings have
seams in Figure \ref{fig:torus2} along the homology basis. Suppose
$z$ and $w$ are two local parameters, differ by a M\"obius
transformation, then $\mu$ should satisfy the following consistency
relation (See Definition \ref{def:discretebeltrami}):
\[
    \mu(w) \frac{\overline{w_z}}{w_z} = \mu(z).
\]
For example, let $p\in S$ is a point on $a_k$ $p \in a_k$, it has two parameters $z_p\in a_k^+$ and $w_p\in a_k^-$, $w_p = \alpha_k(z_p )$. Then $\mu(z_p)$ and $\mu(w_p)$ should satisfy the above
consistent constraint. In our experiments, we find a n-ring neighbor
($n=4$) of $a_k,b_k$, denoted as $R$, then define $\mu(v_i)=z_0$ for
all vertices $v_i$ not in $R$, $\mu(v_i) = 0$, for $v_i$ in $a_k$ or
$b_k$. $\mu$ is extended to other vertices as a complex valued
harmonic function, $\Delta \mu(v_i)=0$, $v_i \not\in R \cup a_k \cup
b_k$, where $\Delta$ is the Laplace-Beltrami operator of the
original surface. This will ensure the consistency relation holds
for $\mu$. In Figure \ref{fig:torus2}, $z_0$ is $0.2+i0.2$ for the left frame of
second row, $0.3$ for the right frame of the second row, and $z_0=z$ for the last row.
From the figure, we can see the deformation of the conformal
structure of the surface (shape of the fundamental domain)with
different Beltrami coefficients. The results show that our method can be applied effectively on general Riemann surfaces of high genus.

\begin{table}
\caption{Computational time.\label{tab:time}}
\begin{center}
\begin{tabular}{ccccccc}
\hline\hline
Figure & \#vertices &\#faces & time(s)\\
\hline
\ref{fig:face4points}&20184& 39984 &99\\
\hline
\ref{fig:face4pointscompose}& 25220 & 49982 & 131\\
\hline
\ref{fig:riemann_mapping}& 20184 & 39984 & 101\\
\hline
\ref{fig:facemultiply1}& 15306 & 29990 & 87\\
\hline
\ref{fig:facemultiply3}& 13515 & 26304 & 108\\
\hline
\ref{fig:kitty}& 10000 & 20000 & 25\\
\hline
\ref{fig:torus2}& 2057 & 4118 & 31\\
\hline\hline
\end{tabular}
\end{center}
\end{table}

\section{Conclusion}
Many surface mappings are quasi-conformal in the real world. According to quasi-conformal Teichm\"uller theory, in general, there exists an one-to-one map between the quasi-conformal maps and the Beltrami coefficients. This work introduces a method to compute quasi-conformal map from Beltrami differentials using auxiliary metric method. The auxiliary metric is constructed from the Beltrami differential, such that the desired quasi-conformal map becomes a conformal one under the auxiliary metric. Theoretic proof shows the rigor of the method, experimental results demonstrate the generality and accuracy of the algorithm. Auxiliary metric method can also be applied to most existing conformal parameterization methods.

\section*{APPENDIX}
\begin{center}
\underline{$\mathbf{PROOF\ OF\ THEOREM\ 4.7 }$}
\end{center}

\bigskip

\begin{proof}
We first compute the entries in the Hessian matrix explicitly. Note that
\[
h_{ij} = \frac{\partial ^2 E }{\partial u_i \partial u_j} = -\frac{\partial K_i}{\partial u_j}
\]
Since $K_i = 2\pi - \sum_{[i,j,k]\in F} \theta_i^{jk}$, differentiating both sides with respect to $u_j$ we get
\[
h_{ij} = -\frac{\partial K_i}{\partial u_j} = \frac{\partial \theta_i^{jk}}{\partial u_j} + \frac{\partial \theta_i^{lj}}{\partial u_j}
\]
Similarly,
\[
h_{ii} = -\frac{\partial K_i}{\partial u_i} = \sum_{j,k}\frac{\partial \theta_i^{jk}}{\partial u_i}
\]
Now, $\frac{\partial \theta_i^{jk}}{\partial u_i}$ and $\frac{\partial \theta_i^{jk}}{\partial u_j}$ can be computed explicitly. We will show the Euclidean case, whereas the hyperbolic case can be proven similarly. Firstly, we observe that by cosine rule
\[
2l_j l_k \cos \theta_i = l_j^2 + l_k^2 - l_i^2
\]
Differentiating both sides with respect to $l_i$,
\[
-2l_j l_k \sin \theta_i \frac{\partial \theta_i}{\partial l_i} = -2l_i
\]
Thus, $\frac{\partial \theta_i}{\partial l_i} = -l_i/2A$, where $A$= area of the triangle $=2l_j l_k \sin \theta_i$
Differentiating both sides with respect to $l_k$,
\[
-2l_j \cos \theta_i - 2l_jl_k \sin \theta_i \frac{\partial \theta_i}{\partial l_k} = 2l_k
\]
Thus, $\frac{\partial \theta_i}{\partial l_k} = \frac{l_j\cos \theta_i - l_k}{2A} = -l_i\cos \theta_j/2A$, since $l_i\cos \theta_j + l_j \theta_i = l_k$.
Also, note that $\frac{\partial l_i}{\partial u_i} = 0$ and $\frac{\partial l_i}{\partial u_j} = l_i$ since $l_i = e^{u_j} e^{u_k} l_i^0$. We have
\begin{equation}\nonumber
\frac{\partial \theta_i}{\partial u_j} = \frac{\partial \theta_i}{\partial l_i}\frac{\partial l_i}{\partial u_j} + \frac{\partial \theta_i}{\partial l_k}\frac{\partial l_k}{\partial u_j} = l_i(l_i - l_k\cos \theta_j)/2A = \cot \theta_k
\end{equation}
Similarly,
\[
\frac{\partial \theta_i}{\partial u_i} = \frac{\partial \theta_i}{\partial l_j} \frac{\partial l_j}{\partial u_i} +\frac{\partial \theta_i}{\partial l_k} \frac{\partial l_k}{\partial u_i} = (-l_i l_j\cos \theta_k + l_i l_k \cos\theta_j)/2A = \cot \theta_j + \cot \theta_k
\]
\end{proof}

\begin{center}
\underline{$\mathbf{PROOF\ OF\ THEOREM\ 4.8 }$}
\end{center}
\begin{proof}
The local convexity of the Euclidean Yamabe energy is due to the fact that the Hessian of $E(u)$ is the matrix $(h_{ij})$ which is semi-positive definite. Clearly, the summation of each row is zero and only the diagonal elements are positive. Furthermore, since the matrix is positive definite on the linear space $\sum_i u_i = 0$, it follows that $H$ is locally strictly convex on the planes. For detail, please see \cite{Luo04}.

\bigskip

We now prove the convexity of the hyperbolic Yamabe energy. This fact was also known to Springborn et al. ahead of us.
We prove the Hessian matrix of the hyperbolic Yamabe energy is positive definite. Let $a_1, a_2, a_3$ be the lengths of a hyperbolic triangle. Make conformal change to produce a new hyperbolic triangle of lengths $y_1, y_2, y_3$ so that $\sinh (y_i/2) = \sinh (a_i/2)e^{u_j + u_k}$, $\{i,j,k\}=\{1,2,3\}$. $\theta_i$ represents the corner angle at the vertex $v_i$. Let $H=[\frac{\partial \theta_i}{\partial u_j}]$ be the matrix.

\bigskip

\noindent $\mathbf{Fact\ 1:}$ $det(H)\neq 0$ for all $u$'s. Indeed, the map from $(u_1,u_2,u_3) \to (\theta_1,\theta_2,\theta_3)$ is a diffeomorphism.

\bigskip

\noindent $\mathbf{Fact\ 2:}$ For any $(a_1,a_2,a_3)$m the set of all $u=(u_1,u_2,u_3)$ such that $ y_i/2 = 2\sinh ^{-1}(\sinh (a_i/2)e^{u_j + u_k})$ satisfies triangular inequalities: $y_1 + y_2 > y_3$, $y_1 + y_3 > y_2$ and $y_2 + y_3 > y_1$ form a connected set $\Omega$ in $\mathbb{R}^3$.

\bigskip

\noindent $\mathbf{Fact\ 3:}$ Since $\Omega$ is connected and $H$ is symmetric in $\Omega$ so that $det(H)\neq 0$, the signature of $H$ is a constant.

\bigskip

\noindent $\mathbf{Fact\ 4:}$ Choose those $u_1, u_2, u_3$ so that $y_1=y_2=y_3$, we see easily by computation that the Hessian $H$ is positive definite.

\bigskip
Thus, $H$ is positive definite over all $\Omega$.

\bigskip
Now to prove $\mathbf{Fact\ 2}:$, introduce a new variable $t_i = e^{u_j + u_k}$. The map $u\mapsto t$ is a diffeomorphism. Thus, it suffices to prove that the set $\Omega^1 = \{(t_1, t_2, t_3)\in \mathbb{R}^3_{>0} | y_i = 2\sinh ^{-1}(t_i\sinh (\frac{a_i}{2}))$ satisfies triangular inequalities $\}$ is connected. Fix $t_2, t_3$, we will show that the set of all $t$ such that
\[
|\sinh^{-1}(t_2\sinh(\frac{a_2}{2})) -  \sinh^{-1}(t_3\sinh(\frac{a_3}{2}))| < \sinh^{-1}(t_2\sinh(\frac{a_2}{2})) + \sinh^{-1}(t_3\sinh(\frac{a_3}{2}))
\]
\noindent is connected. This is obvious since $f(t) = \sinh^{-1}$ is a strictly increasing function in $t$.
\end{proof}

\bigskip

\begin{center}
\underline{$\mathbf{PROOF\ OF\ THEOREM\ 4.9 }$}
\end{center}
\begin{proof}
The solution $u(t) = (u_1(t),..., u_N (t))$ of the discrete Yamabe flow exists for all time so that there are no singularities forming at time equal to infinity. This means that $u_i(t)$'s are in some compact interval in $\mathbb{R}_{>0}$ and also all inner angles $\theta_i^{ij}(t)$ are in some compact interval inside the interval $(0,\pi)$. The matrix $(c_{ij})= (\frac{\partial \theta_i}{\partial u_j})$ has properties that the sum of entries in every row is zero and the diagonal entries are negative. $(c_{ij})$ is symmetric and semi-negative definite. This implies that there is a positive constant $\lambda$ so that the eigenvalues of $(c_{ij})$ considered as a bilinear form restricted to the subspace $\{w\in \mathbb{R}^N| w_1 + ... + w_N = 0\}$ is always bounded by $-\lambda$ for all time $t\in [0,\infty)$, i.e.,
\[
\sum_{i,j} c_{ij} w_i w_j \leq - \lambda \sum_i w_i ^2
\]
when $\sum_{i=1}^N w_i = 0$.

Note that
\[
\begin{split}
\frac{d K_i(t)}{dt} = \frac{d}{dt} (2\pi - \sum_{j,k} \theta_i^{jk}) &= -\sum_{j,k} \frac{d\theta_i^{jk}}{dt}\\
&= -\sum_j \frac{\partial \theta_i}{\partial u_j}\frac{du_j}{dt} = \sum_j \frac{\partial \theta_i}{\partial u_j}(K_j - \overline{K}_j)
\end{split}
\]

Now, consider $G(t) = \sum_{i=1}^N (K_i(t) - \overline{K}_i)^2$. Its derivative can be calculated as
\[
G'(t) = 2\sum_{i,j} c_{ij} (K_i - \overline{K}_i)(K_i - \overline{K}_i)
\]
We have, $G'(t) \leq -\lambda G(t)$. Thus, $G(t) \leq Ce^{-\lambda t}$ and so
\[
|K_i(t) - \overline{K}_i|\leq c_1 e^{-c_2 t}
\]

\end{proof}

\bibliographystyle{IEEEtran}
\bibliography{IEEEabrv,ref}

\end{document}